\documentclass[5p, times]{elsarticle}

\usepackage{amsmath}
\usepackage{xcolor}
\usepackage{bm}
\usepackage{dcolumn,amsthm}
\usepackage{hyperref}
\usepackage{tikz}

\newcommand\restr[2]{{% we make the whole thing an ordinary symbol
  \left.\kern-\nulldelimiterspace % automatically resize the bar with \right
  #1 % the function
  \littletaller % pretend it's a little taller at normal size
  \right|_{#2} % this is the delimiter
  }}

\usepackage{url}
\usepackage{upgreek}
\usepackage{svg}
\usepackage{chngcntr}
\usepackage{apptools}
\newcommand{\littletaller}{\mathchoice{\vphantom{\big|}}{}{}{}}

\AtAppendix{\counterwithin{lemma}{subsection}}

\usepackage{amssymb}
\journal{-}

\usepackage{amsmath,amssymb}

\usepackage{booktabs}

\usepackage{caption}

\usepackage{siunitx}

\usepackage{listings}

\usepackage{url}

\NewDocumentCommand\TODO{O{red}m}{\textcolor{#1}{TODO: #2}}

\begin{document}

\begin{frontmatter}

%% Title, authors and addresses

%% use the tnoteref command within \title for footnotes;
%% use the tnotetext command for theassociated footnote;
%% use the fnref command within \author or \address for footnotes;
%% use the fntext command for theassociated footnote;
%% use the corref command within \author for corresponding author footnotes;
%% use the cortext command for theassociated footnote;
%% use the ead command for the email address,
%% and the form \ead[url] for the home page:
%% \title{Title\tnoteref{label1}}
%% \tnotetext[label1]{}
%% \author{Name\corref{cor1}\fnref{label2}}
%% \ead{email address}
%% \ead[url]{home page}
%% \fntext[label2]{}
%% \cortext[cor1]{}
%% \affiliation{organization={},
%%             addressline={},
%%             city={},
%%             postcode={},
%%             state={},
%%             country={}}
%% \fntext[label3]{}

%\title{Dynamic Effect of Non-Newtonian Cerebral Arterial Circulation on Electrical Conductivity in a Realistic Multi-Compartment Head Model}

\title{A Coupled Diffusion Approximation for Spatiotemporal Hemodynamic Response and Deoxygenated Blood Volume Fraction in Microcirculation}

%% use optional labels to link authors explicitly to addresses:
%% \author[label1,label2]{}
%% \affiliation[label1]{organization={},
%%             addressline={},
%%             city={},
%%             postcode={},
%%             state={},
%%             country={}}
%%
%% \affiliation[label2]{organization={},
%%             addressline={},
%%             city={},
%%             postcode={},
%%             state={},
%%             country={}}
%%%%%%%%%%%%%%%%%%%%%%%%%%%%%%%%%%%%%%%%%
\author[inst1]{Maryam Samavaki\corref{cor1}}
\ead{maryamolsadat.samavaki@tuni.fi}
\cortext[cor1]{Corresponding author at: Sähkötalo building, Korkeakoulunkatu 3, Tampere, 33720, FI}

\affiliation[inst1]{organization={Mathematics, Computing Sciences, Tampere University},%Department and Organization
            addressline={Korkeakoulunkatu 1}, 
            city={Tampere University},
            postcode={33014}, 
            country={Finland}}

            \affiliation[inst2]{organization={Faculty of  Mathematics, K. N. Toosi University of Technology},%Department and Organization
            addressline={Mirdamad Blvd, No. 470}, 
            city={Tehran},
            postcode={1676-53381}, 
            country={Iran}}

\author[inst1]{Santtu Söderholm}

\author[inst2]{Arash Zarrin Nia}

\author[inst1]{Sampsa Pursiainen}
%%%%%%%%%%%%%%%%%%%%%%%%%%%%%%%%%%%%%%

\begin{abstract}
\textbf{Background and Objective:} 
This proof of concept study investigates mathematical modelling of blood flow and oxygen transport in cerebral microcirculation, focusing on understanding hemodynamic responses. By coupling oxygen transport models and blood flow dynamics, the research aims to predict spatiotemporal hemodynamic responses and their impact on blood oxygenation levels, particularly in the context of deoxygenated and total blood volume (DBV and TBV) fractions.
%This article concerns a diffusion-based mathematical model for analyzing blood flow and oxygen transport within the capillaries, emphasizing its significance in understanding the physiological and biochemical dynamics of the cerebrovascular system and brain tissue. The focus of this study is, in particular, on neurovascular coupling and the spatiotemporal aspects of blood flow and oxygen transport in microcirculation.

\noindent \textbf{Methods:} A coupled spatiotemporal model is developed using Fick's law for diffusion, combined with the hemodynamic response function derived from a damped wave equation. The diffusion coefficient in Fick's law is based on Hagen-Poiseuille flow, and arterial blood flow is approximated numerically through pressure--Poisson equation (PPE). The equations are then numerically solved with the finite element method (FEM). Numerical experiments are performed on a high-resolution 7-Tesla Magnetic Resonance Imaging (MRI) dataset for head segmentation, which facilitates the differentiation of arterial blood vessels and various brain tissue  compartments.
%By adopting a coupled spatiotemporal model approach that integrates the hemodynamic response function (HRF) with Fick's law and the Navier--Stokes equations (NSEs), we provide a computational framework for the diffusion-driven transport of deoxygenated and total blood volume fractions (DBV and TBV), essential for understanding blood oxygenation level-dependent functional magnetic resonance imaging (fMRI) and near-infrared spectroscopy (NIRS) applications.

\noindent \textbf{Results:} Numerical analysis reveals how simulated hemodynamic responses influence DBV and TBV fractions, highlighting the complex interactions between blood flow dynamics and oxygen transport. The results demonstrate the model's effectiveness in estimating blood volume fractions, showing its potential in mathematical spatiotemporal modelling of hemodynamic responses, which is in a key role in medical imaging modalities such as blood oxygenation level dependent functional MRI and near infrared spectroscopy. 
%The applicability of the model is further demonstrated through numerical experiments utilizing a 7 Tesla magnetic resonance imaging (MRI) dataset for head segmentation, which facilitates the differentiation of arterial blood vessels and various brain tissue compartments. By simulating hemodynamical responses and analyzing their impact on volumetric DBV and TBV, this study offers valuable insights into spatiotemporal modelling of brain tissue and blood flow.

\noindent \textbf{Conclusions:} This study utilizes spatiotemporal modelling integrated into a realistic head model derived from high-resolution 7 Tesla-MRI to investigate the intricate relationships among cerebral blood flow, oxygen transport, and brain tissue dynamics. This approach not only enhances understanding of cardiovascular conditions but also improves the accuracy of simulations, providing valuable insights into physiological and hemodynamic responses within the human brain. The model's applicability to MRI-based head segmentations suggests potential clinical utility for developing targeted medical interventions.

%The study provides insights into the physiological mechanisms underlying cerebral blood flow and oxygen transport, offering a valuable tool for understanding cardiovascular conditions. By integrating spatiotemporal modelling within a realistic head model derived from high-resolution 7 Tesla-MRI, we analyze the complex interplay between blood flow, oxygen transport, and brain tissue dynamics. The model's applicability to MRI-based head segmentations enhances its potential for clinical use in developing targeted medical interventions. By integrating spatiotemporal modelling within a realistic head model derived from high-resolution 7 Tesla-MRI, we analyze the complex interplay between blood flow, oxygen transport, and brain tissue dynamics. This inclusion of a realistic head model not only enriches the accuracy of our simulations but is also beneficial for understanding the physiological and hemodynamic responses within the human brain.

\end{abstract}

%% Word: 238 words, 1,707 characters

%%Graphical abstract
%\begin{graphicalabstract}
%\includegraphics{grabs}
%\end{graphicalabstract}

%%Research highlights
%\begin{highlights}
%\item Research highlight 1
%\item Research highlight 2
%\end{highlights}

\begin{keyword}
%% keywords here, in the form: keyword \sep keyword
%Blood flow modelling; 
 BOLD effect; cerebral blood flow (CBF);  deoxygenated blood volume (DBV); Fick's law;  hemodynamic response; mathematical model; pressure--Poisson equation (PPE)
%% PACS codes here, in the form: \PACS code \sep code
%\PACS 0000 \sep 1111
%% MSC codes here, in the form: \MSC code \sep code
%% or \MSC[2008] code \sep code (2000 is the default)
%\MSC 65M60 \sep 76D07 \sep 70G45 
\end{keyword}

\end{frontmatter}

%% \linenumbers

%% main text
\section{Introduction}
\label{sec:introduction}
%This study focuses on mathematical modelling of blood flow in capillaries \cite{caro2012mechanics, samavaki2023pressure}, which is essential for a precise understanding of blood and  oxygen transport phenomena in the cerebrovascular system and brain tissue associated with several physiological processes, among other things, gas and nutrition exchange and cell interactions. The transport of oxygen through the microcirculation to the cerebral cells is a crucial element required for maintaining normal brain function and is essential for the biochemical reactions involved in cell metabolism, particularly in converting glucose into adenosine triphosphate (ATP) which is the primary energy source for cellular activities. Unlike muscle tissue, the brain lacks efficient oxygen storage, making it sensitive to fluctuations in oxygen levels. In just a few minutes, an absence of oxygen caused, for example, by a decrease in blood flow can cause the injury and death of brain cells, resulting in neuronal damage. On the one hand, low oxygen levels lead to insufficient energy production; on the other hand, excessively high oxygen levels can lead to oxidative stress, damaging brain tissue \cite{acker2004cellular, nordberg2001reactive}. While only a small amount of oxygen is directly dissolved in the plasma, the majority is carried by hemoglobin, which makes up about $45\%$ of the volume of the red blood cells (RBCs). Hemoglobin plays a significant role in hemodynamics, as it is crucial for transporting oxygen throughout the body. 

In this proof of concept study, we developed the mathematical model that focuses on hemodynamic responses and is divided into two parts of investigation. The first part addresses the spatiotemporal modelling of oxygen transport and blood flow in microcirculation, particularly how cerebral blood flow (CBF) is regulated locally due to neural activity \cite{Brian_2015, Zonta_2003}. The second part explores the numerical analysis of the effect of simulated hemodynamical responses on deoxygenated blood volume (DBV) and total blood volume (TBV) fractions \cite{berg2020modelling, buxton2004modeling}, examining how changes in blood flow dynamics influence these blood volume fractions. A key goal of this study is to address the gap in existing  spatiotemporal hemodynamic response models, which concern simple convolution schemes with constant empirical point spread approximations; in contrast, our study emphasizes mathematical modelling of the diffusion-driven Fick's law  as a parabolic partial differential equation with a variable diffusion coefficient, which previous studies often overlooked  \cite{aquino2014spatiotemporal,lacy2022cortical, Belkhatir_2019_Kalman_filter,pang2016response,pang2017effects,pang2018biophysically,aquino2014deconvolution, drysdale2010spatiotemporal}.

 The brain, which demands a high amount of energy and consumes about 20\% of the body's oxygen, requires a continuous and sufficient supply of oxygenated blood to function properly. Blood flow within the brain is regulated by various anatomical and physiological mechanisms to ensure any disruption in oxygen delivery. Studying blood flow in microcirculation is crucial for understanding how oxygen travels through the cerebrovascular system to reach brain cells, where it is vital for cellular metabolism and energy production.
%Due to the brain's lack of efficient oxygen storage, any interruption in oxygen supply can rapidly lead to brain cell damage or death, while excessive oxygen can cause oxidative stress, highlighting the delicate balance required for proper brain function \cite{acker2004cellular, nordberg2001reactive}. 

Research indicates that during neuronal activation, increased CBF serves a protective role by ensuring sufficient oxygen delivery. Since the early 20th century, studies have observed that blood flow changes in response to increased neural activity \cite{Brian_2015, Zonta_2003}. This specific correlation between neural activity and CBF is referred to as neurovascular coupling \cite{buxton2004modeling, wolf_2002}. Within the cerebral cortex, arterioles and capillaries play crucial roles in this process, ensuring that increased neural activity meets a corresponding rise in blood flow. This enhanced blood flow facilitates the delivery of oxygen and glucose to active brain regions, thereby supporting ongoing neural function \cite{Sava-2014-Large_arteriolar}.

 After neurons become active, the hemodynamic response occurs, involving changes in blood flow and oxygen levels in the brain as a result of neurovascular coupling. As neurons become active, the balance between oxygenated and deoxygenated hemoglobin in the blood shifts. This change can be detected using techniques like near-infrared spectroscopy (NIRS) \cite{ozaki2021near} or  functional magnetic resonance imaging (fMRI) \cite{faro2010bold}, which measures the blood oxygenation level-dependent (BOLD) signal \cite{Tesler_BOLD_2023, kim_BOLD_2007}. These methods are crucial for studying various physiological processes, such as how gases and nutrients are exchanged or how cells interact, by tracking changes in blood flow or oxygenation as indicators of brain activity.

A change in the oxygen level in the blood significantly impacts hemoglobin, a crucial component of red blood cells. The primary hemoglobin role is to be responsible for carrying oxygen throughout the body. When hemoglobin releases its oxygen, it transforms into deoxygenated hemoglobin. By measuring changes in TBV and DBV, we can assess total hemoglobin levels, which correspond to the cerebral blood volume (CBV)  \cite{wolf_2002}. The DBV fraction indicates the percentage of hemoglobin within a given tissue volume that has delivered its oxygen to the surrounding area and becomes deoxygenated hemoglobin.  

The objective of the present work is to develop a comprehensive coupled model that captures the spatiotemporal dynamics of oxygen transport and blood flow in the microcirculation of the brain. By focusing on the hemodynamic response to neural activity, we aim to better understand the mechanisms of neurovascular coupling, with particular emphasis on the local regulation of CBF. Our approach is motivated by the assumption that the relationship between neural activity and blood flow regulation is linear, leading to the formulation of a damped wave equation for predicting the temporal hemodynamic response \cite{aquino2014spatiotemporal,  Aquino_2012, friston2000nonlinear}. This equation is fundamentally concerned with how CBF changes locally in response to neural activity, which is a key mechanism by which the brain regulates its blood supply.

This model is further combined with Fick's law, which governs the diffusion-driven spatiotemporal behavior of DBV and excess TBV fractions under the principle of mass conservation \cite{samavaki2024NSE, samavaki2023pressure}. Findings from \cite{shmuel2007spatio} on the spatiotemporal point-spread function of the blood oxygenation level-dependent magnetic resonance imaging (BOLD MRI) signal at 7T describe how changes in blood oxygenation related to neural activity are distributed across human gray matter over time. In our model, we observe similar spatial extents in excess TBV and DBV fractions, suggesting that this distribution is governed by diffusion, as described by Fick's Law. These dynamics are then measured as the point-spread function in fMRI experiments. 

By applying Fick's law, we can estimate how variations in local blood supply influence the blood flow in microvessels, in particular, DBV and excess TBV fractions. The diffusion coefficient is derived from the laminar Hagen-Poiseuille flow model for cylindrical tubes \cite{samavaki2024NSE, samavaki2023pressure, caro2012mechanics}, and the input boundary conditions for this coupled diffusion model are obtained by approximating the arterial blood flow numerically through the pressure--Poisson equation (PPE) \cite{samavaki2023pressure}. This approach enables the estimation of how variations in local blood supply influence blood flow in microvessels, with the assumption that convection effects can be neglected at a global scale due to fully randomly oriented  microvessels. 
%We concentrate on coupled modelling of spatiotemporal oxygen transport and blood flow in microcirculation, focusing, in particular, on hemodynamic response, in which CBF is regulated locally due to neural activity. We assume that the relationship between that regulation and activity is linear, which motivates predicting the temporal hemodynamic response as a solution for a damped wave equation   \cite{aquino2014spatiotemporal,friston2000nonlinear,friston2003dynamic}. The resulting hemodynamic response function is combined with Fick's law, i.e., a purely diffusion-driven spatiotemporal model \cite{samavaki2023pressure,samavaki2024NSE}, for DBV and TBV fractions along with the principle of mass conservation. By applying Fick's law, we can estimate how variations in local blood supply influence the blood flow in microvessels \cite{epp2020predicting}, in particular, DBV and TBV. This approach is based on the assumption that the densely packed microvessels within the brain tissue are distributed and oriented randomly, suggesting that convection can be omitted on a global scale. The diffusion coefficient is derived from the laminar Hagen-Poiseuille flow model for cylindrical tubes \cite{caro2012mechanics, samavaki2023pressure, samavaki2024NSE}. An input boundary condition for this coupled diffusion model is obtained by approximating the arterial blood flow numerically through Navier--Stokes equations (NSEs) \cite{samavaki2023pressure, samavaki2024NSE, samavaki+tuomela}. 
In the numerical experiments we use a 7 Tesla MRI dataset \cite{svanera2021cerebrum} which allows for obtaining a rough approximation for both arterial blood vessels and different brain tissue compartments \cite{fiederer2016,samavaki2023pressure}.

\section{Methodology} 
%\subsection{Governing equations: macroscopic model}
\label{sec: Theory}

Although in reality, blood is a non-Newtonian fluid with dynamic viscosity  \cite{samavaki2024NSE, saleem_2022, Edward_1969_non-newtonian_blood}, the Newtonian blood flow model is commonly used in studies focusing on oxygen transport and hemodynamic responses \cite{2010_newtonian_hemodynamic}. In this study, we make the simplifying assumption that blood behaves as a Newtonian fluid with constant viscosity, which is valid in macrocirculation, where the vessel diameters are significantly larger than that of a red blood cell  \cite{samavaki2023pressure, Secomb_2017_newtonian_blood}, omitting the possible effect of local changes in the diameter, e.g.\, due to stenosis \cite{kadhim2022numerical}. %Additionally, we assume that the diffusion approximation, specifically the parameter $\delta$ derived from the Hagen-Poiseuille model (see Section \ref{sec:hagen-Poiseuille}), is a rough approximation with limitations. This is due to the lack of MRI data resolution needed to reconstruct microvessels. Consequently, we rely on Fick's model rather than solving the Navier--Stokes equations (NSEs) for microcirculation \cite{samavaki2024NSE}, and a non-Newtonian model is not employed for the microcirculation in this work.

For the aim of this study, we assume that blood transporting oxygen behaves as an incompressible Newtonian fluid \cite{samavaki2023pressure}. This means we describe it using constant viscosity $\mu$ \cite{viscosity_2024, samavaki+tuomela}, pressure $p$, and velocity ${\bf u}$, which are governed by the following PPE \cite{samavaki2023pressure}:
\begin{subequations}
\begin{align}
&\Delta p=\nabla\cdot (\rho \,{\bf {\bf f}})+2\mu \,\nabla\cdot ({\bf Ri}({\bf u})) &\mathsf{in}\,\Omega \! \times \! [0, T]\,,
\label{eqn:line-1.1}
 \\
& \mathsf{div}({\bf u})=0&\mathsf{in}\,\Omega \! \times \!  [0, T]\,,
\label{eqn:line-1.2}
  \\
&{\bf u}({\bf x};0)=0 &\mathsf{on}\,\Omega \,, 
\\
&{\bf g}(\nabla p, {\vec{{\bf n}}})=-\zeta \lambda  ( p -    {p}^{\mathcal{(B)}})&\mathsf{on}\, \mathrm{B}\,.
\label{robin}
\end{align}
\label{main-SE}
\end{subequations}
In this study, we consider a 3D compact domain $\Omega\subset\mathbb{R}^3$, representing the human brain, and divide it into two distinct regions:
\[
 \Omega =  \Omega_{\mathrm{A}} \cup \hat{\Omega}  \quad \hbox{with} \quad \partial\Omega_{\mathrm{A}}\cap \partial\hat{\Omega} = \mathrm{B}
\]
In this context, $\Omega_{\mathrm{A}}$ represents the arterial blood flow, while $\hat{\Omega}$ denotes the microcirculation domain. The boundary between these two regions, denoted by $\mathrm{B}$, is where blood flows out of the arteries and enters the microcirculation domain. This transition happens in arterioles, which are small blood vessels with relatively thick walls. As blood moves through these arterioles, the pressure drops according to the formula $-{\bf g}(\nabla p, {\vec{{\bf n}}})=-\partial p/\partial\vec{{\bf n}}$, before it continues into the thin-walled capillaries. In order to determine blood pressure $p$ and velocity ${\bf u}$ of incompressible blood flow inside the arteries ${\Omega}_\mathrm{A}$, we use an approximation of the solution of regularized PPE as shown in \cite{samavaki2023pressure}. In the microcirculation region $\hat{\Omega}$, blood flows through capillaries before it moves into the veins via venules. Because the diameter of the microvessels, including arterioles, capillaries, and venules, is smaller than the finest details of what is possible to model numerically, we consider the effect of the arterioles as a boundary condition set on $\mathrm{B}$. Within $\hat{\Omega}$, the blood flow simulates diffusion through porous medium, where the randomly oriented capillaries constitute the pores. 

In this study, another simplifying assumption is that the diffusion approximation, especially the Hagen-Poisseuille model-based diffusion parameter $\delta$, see Section \ref{sec:hagen-Poiseuille}, is a rough approximation since we do not have a capability of reconstructing the microvessels from MRI data and, therefore, we rely on Fick's law rather than solve the Navier--Stokes equations (NSEs) for microcirculation; thus, the non-Newtonian model is not used for microcirculation either.

%%%%%%%%%%%%%%%%%%%%%%
%%%%%%%%%%%%%%%%%%%%%%
%%%%%%%%%%%%%%%%%%%%%%
\section{Modelling of Coupled Blood and Oxygen Diffusion}
\label{sec:Fick_law} 

Brain tissue is composed of various components, including cells, extracellular matrix, and capillary vessels. These components together form a complex structure that behaves like a porous medium for the blood and oxygen \cite{Swartz_2007}. Hemoglobin is a protein in red blood cells responsible for carrying oxygen from the lungs to the tissues and returning carbon dioxide from the tissues to the lungs. Oxygen carried by oxygenated blood  is utilized in the tissues, leading to the formation of DBV, which is an indicator of oxygenation levels in the blood. This DBV then travels through the network represented by $\hat{\Omega}$, enters the venous system through venules, and transports carbon dioxide from the tissues back to the lungs.
%the capillary wall does not consume oxygen and has a lower diffusion coefficient. %We aim to explore how concentration gradients affect oxygen uptake within tissue, their transport within blood vessels, and the movement of mass across BBB. 

 In the capillary bed, oxygen molecules move from the blood to the tissue through the blood-brain barrier, which acts as a selective membrane that controls the passage of substances into and out of the brain \cite{ Kaplan_2020_BBB, Simpson_2006_BBB}. This process can be enhanced locally due to neural activity, which will increase the blood flow supplied to the active region. This local regulation of the microcirculation is known as the hemodynamic response \cite{aquino2014spatiotemporal, friston2000nonlinear, friston2003dynamic}. 

In this section, we investigate the dynamic time evolution of DBV and TBV fractions, represented by the variables $\tilde{\bf{q}}$ and $\tilde{\bf{c}}$, respectively. Of these, TBV can be decomposed into dynamic excess ${\bf c}$ and static background $\overline{\bf c}$ distribution. These variables are approximated by distributions within the domain $\hat{\Omega}$ and are treated as coupled phenomena influenced by neural activity through the hemodynamic response.
%In this study, we assume DBV represented by $\tilde{\bf{q}}$. The DBV is related to $\tilde{\bf{c}}$ through a factor $q, \tilde{\bf{q}}=q \tilde{\bf{c}}$, where $q$ represents the local DBV fraction, relative to the TBV fraction $\tilde{\bf{c}}=(\tilde{c}, \tilde{c}, \tilde{c})$. Here $\tilde{\bf{c}}=\bf{c}+\overline{\bf{c}}$ where $\bf{c}$ stands for the regular volume fraction and $\overline{\bf{c}}$ represents a constant background concentration. Both $\bf{c}$ and $\overline{\bf{c}}$ are also 3D vectors; $\bf{c}=(c, c, c)$ and $\overline{\bf{c}}=(\bar{c}, \bar{c}, \bar{c})$, respectively. All these values are averaged over a small volume.  Once oxygen molecules reach the tissue, they further diffuse into the tissue matrix. However, there is typically no detection term associated with this process in Fick's Law \cite{}, as there is no bulk movement of oxygen within the tissue itself. Instead, diffusion is the dominant mechanism driving oxygen movement within the tissue in the sense that oxygen molecules move from areas of higher concentration to areas of lower concentration due to random molecular motion, which is the essence of diffusion.
%\subsection{Convection-Diffusion Equation}
To model the complex process of oxygen transport and consumption within the porous domain $\hat{\Omega}$, which includes both the microvessel network and tissue, we employ the convection-diffusion reaction equation to describe the DBV fraction, $\tilde{{\bf q}}$. The equation is formulated as follows \cite{aquino2014spatiotemporal, lacy2022cortical,pang2016response,pang2017effects}: 
\begin{equation}
\begin{aligned}
& \tilde{{\bf q}}_{,t}=-\mathsf{div}({\bf u}\otimes \tilde{{\bf q}})+\mathsf{div}({\bf J}_{\tilde{{ q}}}) + {\bf s} & \mathsf{in}\,\,\hat{\Omega} \! \times \! [0, T]\,,
 \\
& \nabla\cdot{\bf u}=0 & \mathsf{in}\,\,\hat{\Omega} \! \times \! [0, T]\,.
\end{aligned}
\label{DH_concentration}
\end{equation}
In this system of equations, the first term on the right-hand side of the first equation represents advection, which describes the transport of deoxygenated blood by the velocity field ${\bf u}$. The second term accounts for diffusion, determined by the flux density ${\bf J}_{\tilde{ q}} = - \delta \, \nabla \tilde{\bf q}$, with $\delta$ denoting a scale-dependent effective diffusion coefficient (see Section \ref{sec:hagen-Poiseuille}). The third term, denoted by $\bf s$,  corresponds to the sink/source contributions, which describe oxygen consumption by the brain tissue and the outflow of blood to the venous system. The second equation enforces the incompressibility condition, ensuring that the blood flow is divergence-free throughout the domain. This implies that the volume of blood in any region remains constant over time, behaving as an incompressible fluid. Consequently, the amount of blood entering any region exactly matches the amount leaving, with no accumulation or loss.

% Thus
% \begin{subequations}
% \begin{align}
% & {\bf c}_{\mathrm{b},t}+{\bf u} \cdot \nabla{\bf c}_{\mathrm{b}}-\nabla \cdot(D_{\mathrm{b}} \, \nabla {\bf c}_{\mathrm{b}})={\bf u} \cdot \nabla{\bf c}_{\mathrm{T}}-\nabla \cdot(D_{\mathrm{T}} \, \nabla {\bf c}_{\mathrm{T}})
% \\
% & \nabla\cdot{\bf u}=0
% \end{align}
% \label{blood_o2}
% \end{subequations}
% where ${\bf c}_{\mathrm{b}}$ states as a blood oxygen concentration, $?$ as the intensity of oxygen exchange between the blood and tissue, and $\bf u$ the velocity field of the capillary network.
%%%%%%%%%%%%%%%%%%%%%%%%%%%%
%%%%%%%%%%%%%%%%%%%%%%%%%%%%
%%%%%%%%%%%%%%%%%%%%%%%%%%%%

\subsection{Background Blood Volume Fraction}

We first derive an approximation for the background TBV, $\overline{\bf c}$ based on the length density microvessels $\xi$  ($m^{-2}$) \cite{kubikova2018numerical}, i.e., the number of microvessels per unit area, which varies between different tissue types and is the sum $\xi = \xi_a + \xi_c + \xi_v$ of the length densities of arterioles, capillaries, and venules $\xi_a$, $\xi_c$, and $\xi_v$, respectively (Figure \ref{fig:artery_arteriole}). To find this sum, we consider the corresponding  cross-sectional areas $A_a$, $A_c$ and $A_v$ of each vessel type, the total cross-sectional area of all the microvessels $A$, and area fractions $\gamma_a$, $\gamma_c$, and $\gamma_v$, $\gamma_a + \gamma_c + \gamma_v = 1$ which are linked as follows:
\[
\frac{A_a}{\gamma_a} \xi_a = \frac{A_c}{\gamma_c} \xi_c = \frac{A_v}{\gamma_v} \xi_v = A.
\]

As a result, we find that the arteriolar length density can be determined by taking the inverse of the sum of the cross-sectional area ratios and their respective area fractions:
\[
\xi_a = \xi \left( 1 + \frac{A_a \gamma_c}{A_c \gamma_a} + \frac{A_a \gamma_v}{A_v \gamma_a} \right)^{-1}\,.
\]
From this expression, we can deduce that the densities of capillaries $\xi_c$ and venules $\xi_v$ are related to the density of arterioles $\xi_a$ by the ratios of their cross-sectional areas and area fractions as follows:
\[
\xi_c = \frac{A_a \xi_a \gamma_c}{A_c \gamma_a}\,, \quad \xi_v = \frac{A_a \xi_a \gamma_v}{A_v \gamma_a}\,.
\]
Using the above expressions, background distribution is modelled as 
\[
\overline{\bf c} = \theta \,  (\xi_a A_a + \xi_c A_c + \xi_v A_v), 
\]
i.e., it is the total approximate volume fraction of microvessels multiplied by $\theta$, which denotes the relative expansion of microvessels' cross-sectional area due to vessel upregulation.

%%%%%%%%%%%%%%%%%%%%%%%%%%%
%%%%%%%%%%%%%%%%%%%%%%%%%%%
%%%%%%%%%%%%%%%%%%%%%%%%%%%
\subsection{Hagen-Poiseuille Model for Arterioles}
\label{sec:hagen-Poiseuille}

\begin{figure}[h!]
\centering
\includegraphics[height=3cm]{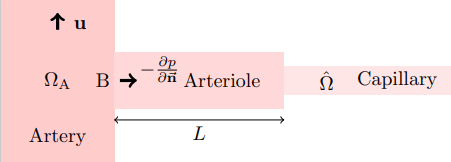}
%\begin{figure}[ht!]
%\centering
%\input{artery-arteriole-image.tex}
\caption{Computational domain in cylindrical coordinates, taking into
account for the axial symmetry in which blood flow in the artery domain $\Omega_\mathrm{A}$ is determined by the velocity field ${\bf u}$. Blood flowing through the boundary $\mathrm{B}$ is aligned along with normal vector $\vec{\bf n}$ and enters first arterioles with length $L$, where the pressure $p$ drops linearly, i.e., $\partial p / \partial \vec{\bf n}$ is constant, before it enters the capillary bed $\hat{\Omega}$.}
\label{fig:artery_arteriole}
\end{figure}

Figure~\ref{fig:artery_arteriole} presents the passage of blood from artery $\Omega_A$ into the capillaries $\hat\Omega$ via arterioles. Blood crosses the boundary $B$ as it moves from the arteries into the arterioles, which serve to reduce the blood pressure before it reaches the capillaries.

%%%%%%%%%%%%%
%%%%%%%%%%%%%
%%%%%%%%%%%%%
%%%%%%%%%%%%%%
%%%%%%%%%%%%%%
%%%%%%%%%%%%%%
  
%\textcolor{black}{
In the core of our approach is the Hagen-Poiseuille model  \cite{caro2012mechanics}, which we use to approximate the transition of blood flow  from macro to microcirculation via arterioles and, consequently, to derive a uniform approximative diffusion coefficient $\delta$ for microcirculation \cite{samavaki2024NSE, samavaki2023pressure}. We consider each arteriole as a cylindrical pipe of length $L$, where the flow is laminar, and we assume that the flow parameters are either constant or vary linearly along its length (Figure \ref{fig:artery_arteriole}). Specifically, we assume that both $\nabla \tilde{q}$ and $\nabla \tilde{c}$ are constant, with proportional relationship $\nabla \tilde{q} = \tau \nabla \tilde{c}$, where $\tau$ is an (implicit) coupling constant and $\tilde{c}$ represents TBV. Given that the pressure drop and TBV are assumed to be linearly proportional to each other, the DBV gradient can be approximated as 
$
|\nabla \tilde{c}| = |\vartheta \tilde{c}/L|
$. Consequently, the boundary derivative of the influx of deoxygenated blood entering from arteries through arterioles can be approximated as $|\mathbf{J}_{ \tilde{{ q}}}| =|-\delta \nabla \tilde{{ q}}|= |-\delta \tau \nabla \tilde{{ c}}|$. This simplifies to:
\begin{equation}
|\mathbf{J}_{ \tilde{ q}}| =  |-\delta  \nabla \tilde{q}| =  \delta \tau \frac{ \vartheta  \tilde{ c}}{L}= \delta \tau \frac{ \vartheta \,\theta \,A_a \,\bar{\xi}_a}{L}\,.
\label{hemoglobin_conservation}
\end{equation}
In this context, the term $A_a \overline{\xi}_a$ represents the mean volume concentration of the blood, $\delta$ is the diffusion coefficient, and $\vartheta$ denotes the relative pressure drop along the arterioles. Additionally, the terms $A_a$ and $L$ correspond to the cross-sectional area and length of a single arteriole, respectively, while $\bar{\xi}_a$ refers to the average length density ${\xi}_a$ ($m^{-2}$) of arterioles per tissue volume. The average length density is given by:
\begin{align*} 
\overline{\xi}_a = \frac{1}{|\partial \Omega_A|} \int_{\partial \Omega_A} \xi_a \,\hbox{d} \omega_{\partial \Omega_A}\,,
\end{align*} 
where $|\partial \Omega_A|=\int_{\partial \Omega}\,\mathrm{d}\omega_{\partial \Omega_A}$.

In this context, $|{\bf J}_{\tilde{q}}| = \tau  \theta  Q_a {\xi}_a$, where the term $Q_a$ represents the flow rate of a single arteriole \cite{samavaki2023pressure}. To approximate $Q_a$, we use the total cerebral blood flow rate $Q$, which is related to $Q_a$ through the relationship $Q = |B|  Q_a  \bar{\xi}_a$. Utilizing the Hagen-Poiseuille equation, we derive an approximation for the length of the arteriole as follows:
\begin{align*}
\vartheta \bar{p} =\frac{8 \pi \mu L Q_a}{A_a^2}\,, \quad \hbox{ i.e., } \quad L = \frac{\vartheta\bar{p} A_a^2}{8 \pi \mu Q_a}\,,
\end{align*}
where $\mu$ is blood viscosity and $\bar{p}$ denotes average pressure. Here, $\vartheta \bar{p}$ represents the pressure drop along the length of the arteriole. Solving for $\delta$ from equation \eqref{hemoglobin_conservation} leads to the formula 
\begin{equation}
\delta  =  \frac{\lambda A_a \bar{p}}{8 \pi \mu} \,,
\label{diffusivity}
\end{equation}
where $\lambda = \xi_a / \overline{\xi}_a$ is the ratio between the length density $\xi_a$ and its average value $\overline{\xi}_a$.
%%%%%%%%%%%%%%%%%%%%%%
%%%%%%%%%%%%%%%%%%%%%%%
%%%%%%%%%%%%%%%%%%%%%%%%
\subsection{Coupled Fick's Law for Deoxygenated Blood Volume Fraction}
\label{Coupled_Fick}
%We assume the velocity field $\mathbf{u}$ and concentration of DBV,  $\tilde{{ \bf q}}$  {\color{red} $\tilde{{\bf c}}_{\mathsf{DBV}}=q\,\tilde{{ \bf c}}$ in which $q$ states as a local DBV mass fraction,} relative to total volumetric blood concentration $\tilde{\bf c} =(\tilde{ c}, \tilde{ c}, \tilde{ c})$ in which $\tilde{\bf c} = {\bf c} + \bar{\bf c}$, where ${\bf c}=(c, c, c)$ and $\bar{\bf c}=(\bar{ c}, \bar{ c}, \bar{ c})$ denote excess and constant background  concentration, respectively,  are averaged over a small volume.

%We assume there is a velocity field represented by $\bf u$ and a concentration of DBV, represented by $\tilde{{ q}}$. 
%DBV is related to TBV $\tilde{{c}}={ c}+\overline{{ c}}$, where ${ c}$ stands for the regular concentration and $\overline{{ c}}$ represents a constant background concentration. %Both ${ c}$ and $\overline{{ c}}$ are also three-dimensional vectors; ${ c}=(c, c, c)$ and $\overline{{ c}}=(\bar{c}, \bar{c}, \bar{c})$, respectively. All these values are averaged over a small volume.

Due to the averaging process, the flow of DBV can be described by the advection-diffusion equation, derived from Fick's second law of diffusion \cite{berg2020modelling, arciero2017mathematical, reichold2009vascular}. This equation characterizes the flow as primarily diffusive, aligning with the principles of Fick's law. According to this law, the flux density of DBV is given by ${\bf J}_{\tilde{q}} = - \delta \nabla \tilde{q}$, where $\delta$ is the diffusion constant as discussed in Section \ref{sec:hagen-Poiseuille}. 

In a single capillary vessel, the blood flow is fully advective, meaning that the movement of blood is driven by the flow itself. In this scenario, diffusion's contribution, which is the process by which substances spread out due to concentration gradients, is minimal, almost negligible. However, when we expand our focus to a small volume of tissue containing many microvessels, these vessels are typically oriented in various random directions. As a result, the velocities of the blood flow in different directions tend to cancel each other out, i.e., $\mathbf{u} = \mathbf{0}$. This means that, on average, there is no net movement of substances in any specific direction due to the flow. In this case, advection does not contribute to the movement of substances at all, so the advection term in the equation becomes zero. Thus, in the system of equations \eqref{DH_concentration}, the term $\mathsf{div}({\bf u}\otimes \tilde{{\bf q}})=0$ means that there is no net divergence in the flux of $\tilde{{\bf q}}$ in the direction of the velocity field ${\bf u}$. As a result, we assume that DBV conservation is completely driven by diffusion \cite{drysdale2010spatiotemporal} and, thereby obeys a conservation equation of the form:
%Namely, while in a single capillary vessel the flow is fully advective with a vanishing diffusion, for a small tissue volume containing fully randomly oriented  microvessels, the averaged velocity and, thereby, the advection term vanishes, i.e., $\mathbf{u} = \mathbf{0}$ and $\mathsf{div}({\bf u}\otimes \tilde{{\bf q}})=0$. Thus, we assume that DBV conservation is completely driven by diffusion \cite{drysdale2010spatiotemporal} and, thereby, obeys a conservation equation of the form:
 \begin{equation}
\begin{aligned}
 &\tilde{{  q}}_{,t} - \mathsf{div}\big(\delta\, \nabla \tilde{{  q}}\big)= -\varepsilon\tilde{{  q}} + ({{ h}}\tilde{c} - \tilde{{  q}} ) \upsilon 
 &\mathsf{in}\,\,\hat{\Omega} \! \times \! [0, T]\,,
 %\label{ficks_law_q_1}
 \\
 &{\tilde { q}}(\cdot  ; 0)= (1-h) {\bar { c}}&\mathsf{in}\,\,\hat{\Omega}\,,
 \\
 &\mathbf{g}( \delta  \nabla {\tilde { q}} ,{\vec{ \bf{n}}}) = -\mathbf{g}( {\bf J}_{\tilde{q}} ,{\vec{ \bf{n}}})= (1-h)f  &\mathsf{on}\,\,\mathrm{B} \,.
 %\label{DE.vessel_2}
% \\
 %&\mathbf{g}(\nabla q, {\vec{ \bf{n}}}) = 0  &\mathsf{on}\,\,\mathrm{B} \,.
 %\label{DE.vessel_3}
 \end{aligned}
 \label{ficks_law_q}
 \end{equation}
%%%%%%%%%%%%%%%%%%
%Here, the second term on the left-hand-side introduces diffusion, with  $\delta$ denoting a scale-dependent effective diffusion coefficient. 
The term on the right-hand side, $s=- \varepsilon \tilde {q} + ({{ h}} \tilde{c}  - \tilde{{  q}})  \upsilon $, represents the sink/source. In this expression, $({{ h}} \tilde{c}  - \tilde{{  q}})  \upsilon$ corresponds to the effect of the oxygen consumption rate. The term $h \tilde{c}$ represents the oxygenated blood volume (OBV) fraction in microcirculation, and here serves as a coupling between the models of TBV and DBV represented by $\tilde{c}$ and $\tilde{q}$ along with the boundary condition. The parameter $h$ denotes the reference OBV fraction in blood \cite{aquino2014spatiotemporal}. The fractional oxygen consumption rate, denoted by $\upsilon = \eta \rho^{-1} \psi^{-1}$, determines the rate at which oxygen is consumed as it passes from oxygenated blood to cortical tissue. This rate is calculated by dividing the molecular hemoglobin consumption rate $\eta$ by both the molarity of hemoglobin in blood $\psi$ and the mass density of blood $\rho$. The term $\varepsilon \tilde q$ denotes the rate at which DBV decreases due to blood outflow. The parameter $\varepsilon$ has been selected so that the sink, which corresponds to the volume $\mathcal{V}_{\max} = (4/3) \pi R^3$ of the largest element (radius $R$) in the finite element, matches the flux  $|\mathbf{J}|  =\delta\tau \vartheta {\tilde{c}} /L $ on its surface. We approximate the volume $\mathcal{V}_{\max}$ of the element as that of a sphere with radius $R$ \cite{samavaki2023pressure}, leading to the following expression:
\begin{equation*}
\varepsilon \mathcal{V}_{\max }=(4 / 3) \pi R^3 \varepsilon=4 \pi R^2|\mathbf{J}|  \quad \text { or } \quad \varepsilon=\delta\tau \frac{\vartheta {\tilde{c}}}{L}\left(\frac{36 \pi}{\mathcal{V}_{\max }}\right)^{1 / 3}.
%\label{varepsilon}
\end{equation*}
%The term $f$ represents the pressure from PPE that causes this inflow. 
The initial distribution $(1-h) {\bar { c}}$ stands as a background DBV fraction. The boundary condition shows the amount of blood flowing into the system. The term $f$ represents the influx of blood into the domain from the arterial region $\Omega_{A}$,  and in our model is derived from the arterial flow described \cite{samavaki2023pressure}. Specifically, $f$ measures the rate at which DBV passes crosses the boundary $B$ and enters the microcirculation domain. Thus, in the boundary condition, $f$ represents the pressure difference which is coming from PPE that causes this inflow, which is expressed as $(1-h)f$, where $(1-h)$ adjusts the total blood flux to account for the fraction of DBV. Mathematically, this is represented by $\mathbf{g}(\mathbf{J}_{\tilde{q}}, {\vec{ \bf{n}}})$, with $\mathbf{J}_{\tilde{q}}$ denoting the flux density of deoxygenated hemoglobin and ${\vec{ \bf{n}}}$ being the outward normal unit vector on the boundary $B$.

Note that deoxygenation from tissue may not solely occur through passive diffusion driven by the oxygen gradient between blood vessels and surrounding tissue. Instead, it might involve an active process. In this model, we integrate an active diffusion mechanism in accordance with Fick's law.

%%%%%%%%%%%%%%%%
%%%%%%%%%%%%%%%%%
%%%%%%%%%%%%%%%%%%
%%%Fick's Law for Microcirculation%%%%%%
%%%%%%%%%%%%%%%%%
\subsection{Fick's Law for Microcirculation}

Fick's law, describing the blood volume fraction of excess TBV, denoted by $c$, compared to the resting state, can be written as:
\begin{equation}
\begin{aligned}
&c_{, t} - \mathsf{div}(\delta \,\nabla c)=-\varepsilon c+\alpha_{, t} & \mathsf{in}\,\, \hat{\Omega} \times[0, T] \,,
%\label{ficks_law_for_microcirculation_1}
\\
&\mathbf{g}(\delta \, \nabla c, {\vec{ \bf{n}}})=- \mathbf{g}(\mathbf{J}_{c}, {\vec{ \bf{n}}})  = f &\mathsf{on}\,\,\mathrm{B}\,. 
%\label{ficks_law_for_microcirculation_2}
\end{aligned}
\label{ficks_law_for_microcirculation}
\end{equation}
The method used to calculate the diffusion coefficient $\delta$ for $c$ is the same as that used to derive $\tilde q$ in Section \ref{sec:hagen-Poiseuille}. Given that the minor influence of blood oxygenation level on blood viscosity is negligible \cite{valant2016influence}, we apply the same diffusion coefficient $\delta$ approximation for both excess TBV and DBV fractions. The rate of blood outflow from capillaries to venules and venous vessels is represented by the coefficient $\varepsilon$ in this formulation, which is computed in the same way as $\tilde q$ in Section \ref{Coupled_Fick}. The term $\alpha$ represents a local regulation term, corresponding to the hemodynamic response, which encompasses the local rate of blood volume inflow and its response to neural activity, while $\mathbf{J}_{c} = -\delta\, \nabla { c}$ denotes the flux density of excess TBV. 

%$\mathbf{J}_{c} = -\delta\, \nabla { c}$ is the flux density of excess TBV; the term  $\varepsilon$ is a coefficient representing the rate of blood outflow from capillaries to venules and venous vessels \cite{samavaki2023pressure}; and $\alpha$ is a local regulation term corresponding to the hemodynamic response, the local rate of blood volume inflow, and the response to neural activity. The diffusion constant, $\delta$, for $c$ can be derived from Fick's law as for $\tilde{q}$ in Section \ref{sec:hagen-Poiseuille}. Omitting the minor effect of blood oxygenation level on blood viscosity \cite{valant2016influence}, we apply the same diffusion coefficient $\delta$ approximation for both TBV and DBV. 

%%%%%%%%%%%%%%%%%%%%%%%
%%%%%%Balloon model%%%%
%%%%%%%%%%%%%%%%%%%%%%

\subsubsection{Balloon model}
\label{sec:balloon_model}

The relationship between changes in blood flow and neuronal activity in the brain is explained by the balloon model \cite{friston2000nonlinear, Buxton_1998}. The hemodynamic response, or a change in blood flow and oxygenation in response to brain activity, is estimated with the use of this model. To study how the local regulation $\alpha$, assumed to vary smoothly over time, affects the local temporal dynamics, we use a modified version of the balloon model. Here, $\alpha$ is represented as $\alpha = \beta r$, where $\beta$ is a mollifier (smoothing) function and $r$ is a hemodynamic response, influenced linearly by the local strength of the neural activity. This approach is based on the damped harmonic oscillator model of the hemodynamic response \cite{aquino2014spatiotemporal, friston2000nonlinear}.
%To estimate the local temporal effect of the local regulation $\alpha$, which we assume to be smooth function of time, we apply a time-mollified version of the so-called balloon model: $\alpha = \beta r$, where $\beta$ is a mollifier function and $r$ is a hemodynamic response that depends linearly on  the  local strength  of the neural activity  \cite{aquino2014spatiotemporal, friston2000nonlinear} leading to the following damped harmonic oscillator model of the hemodynamic response:
\begin{equation}
{ r}_{,tt} +  \kappa { r}_{,t} + \gamma { r} = { \zeta}\,. 
\label{DHOM}
\end{equation}
The term $\kappa { r}_{,t}$ corresponds to the damping force that reduces the oscillatory changes in the blood flow in response to neuronal activity, and $\kappa$ represents the rate at which the flow signal decreases after neural activity. The term $\gamma { r}$ represents the restoring force that acts to return the blood flow to its equilibrium state after a neuronal activity, and $\gamma$ is the flow-dependent elimination constant, which indicates the rate of removing substances from tissue in response to a change in blood flow. In the right-hand side of the equation \eqref{DHOM}, the term ${ \zeta}$ measures oxygen consumption due to neural activity. It acts as an external driving force for the hemodynamic response, influencing the response $r$. The mollifier is defined as 
\[
\beta(t) = 
\exp\left(1 - \frac{1}{1 - \left( \frac{2t - T}{T} \right)^2}\right) 
\]
for the interval $[0, T]$. We introduce a smoothening function, $\beta (t)$, to maintain the smoothness of the hemodynamic response at the endpoints. Without mollification, the derivative at the beginning would display a discontinuity \cite{friston2000nonlinear}, but by applying the mollifier, we provide a smooth derivative, matching the solution with the expected behavior of the real hemodynamic response.

\section{Variational Formulation}
\label{sec:Variational Formulation}

In this section, we address the problem of the flow of an incompressible fluid within a compliant domain. The model of oxygen transport is under consideration in a 3D circular cylindrical tube of blood vessels.  DBV in cerebral arteries is approximated numerically using finite element (FE) discretization. To proceed with FE discretization, we must identify an appropriate variational formulation of Fick's law. Both the mathematical analysis and the numerical solution are based on weak formulations.

We assume that $\tilde{q} \in \mathbb{V}$ and $c \in \mathbb{V}$ with $\mathbb{V}=\mathrm{H}^1(\Omega)$. Here, $\mathrm{H}^1(\Omega)$ is known as the \textit{Sobolev space}, consisting of functions defined on $\Omega$ that are square integrable, i.e., $\int_{\Omega} |{u}|^2 \mathrm{~d}\omega_{\Omega} < \infty$, and whose partial derivatives are also square integrable. Formally, it can be expressed as:
\[
\mathrm{H}^1(\Omega)=\left\{u\in \mathrm{L}^2(\Omega)\,|\, \nabla u\in \mathrm{L}^2(\Omega )\right\}\,.
\]
%%%%%%%%%%%%%%%%%%%
%%%%%%%%%%%%%%%%%%%
%%%%%%%%%%%%%%%%%%%
\subsection{Variational Formulation for Deoxygenated Blood Volume Fraction}
\label{Variational_Formulation_DBV}
Considering blood as a Newtonian fluid \cite{samavaki2023pressure}, the diffusion coefficient $\delta$ remains constant, as defined by equation \eqref{diffusivity}. Now, let's express the variational formulation of the system described by (\ref{ficks_law_q}):
\begin{quote}
    Find ${\tilde q}\in {\mathbb{V}}$ such that, for a smooth enough test function, $ \varphi\in  {\mathbb{V}}$
    \begin{equation}
    \begin{aligned}
    d({\tilde q}, \varphi)
    =
    \int_{ \mathrm{B}}
    (1-h) f \varphi
    \,\mathrm{d} \omega_{\mathrm{B}}
    +
    \int_{\hat\Omega}
    h\upsilon\tilde{c}\ \varphi
    \,\mathrm{d}\omega_{\hat\Omega}
    \,.
    \end{aligned}
    \label{q_fick_variational_2}
    \end{equation}
\end{quote}
In this expression, the term representing the linear boundary for the incoming flow is located on the left-hand side. % while $\vec{\bf n}$ stands for the normal unit vector defined within the microcirculation domain. 
The continuous bilinear form $d(\cdot, \cdot)$ is defined as $d:{{\mathbb{V}}}\times {\mathbb{V}}\rightarrow \mathbb{R}$, where
\begin{equation}
\begin{aligned}
d(\tilde{q}, \varphi) &= \frac{\partial}{\partial t}\int_{\hat\Omega}\tilde{q} \,\varphi\,\mathrm{d}\omega_{\hat\Omega}+\int_{\hat\Omega} \delta\,{\bf g}(\nabla \tilde{q}, \nabla \varphi)\,\mathrm{d}\omega_{\hat\Omega}
\\
&+\int_{\hat\Omega} \varepsilon \tilde{q}\varphi\,\mathrm{d}\omega_{\hat\Omega} + \int_{\hat\Omega} \tilde{q}\upsilon\varphi\,\mathrm{d}\omega_{\hat\Omega}\,.
%\\
%&-\int_{\hat\Omega} \alpha\zeta  \varphi\,\mathrm{d}\omega_{\hat\Omega}\,.
\end{aligned}
\label{q_fick_variational_g}
\end{equation}
\subsection{Variational Formulation for Excess Total Blood Volume Fraction}
\label{Variational_Formulation_TBV}
By applying equation \eqref{ficks_law_for_microcirculation} to a chosen test function $\varphi$ from the subspace $\mathbb{V}$, and subsequently using the method of integration by parts, we obtain the variational form of Fick's law for the volume fraction of excess TBV, denoted by $c$. The resulting expression is as follows:
\begin{quote}
  Find $c\in {\mathbb{V}}$ such that, for a smooth enough  test function, $ \varphi\in  {\mathbb{V}}$ 
\begin{equation}
    a(c,\varphi)
    =
    \int_{\mathrm{B}}
    f \varphi
    \, \mathrm{d} \omega_{\mathrm{B}}
    +
    \int_{\hat\Omega}
    \alpha_{, t} \,\varphi
    \,\mathrm{d}\omega_{\hat\Omega}
    \,.
    \label{variation_2_part_1}
\end{equation}
\end{quote}
Within this context, \( \vec{\mathbf{n}} \) symbolizes the normal unit vector on boundary of the microcirculation area, with the incoming flux density's linear term presented on the left side. The bilinear form $ a(\cdot,\cdot)$ is defined through mapping $a: \mathbb{V} \times \mathbb{V} \rightarrow \mathbb{R} $, detailed as
{\setlength\arraycolsep{2 pt} \begin{eqnarray}
    a(c,\varphi) & = & \frac{\partial}{\partial t} \int_{\hat{\Omega}} c \varphi \,\mathrm{d}\omega_{\Omega} + \int_{\hat\Omega} \delta g( \nabla c, \nabla \varphi) \,\mathrm{d}\omega_{\hat\Omega} \nonumber \\ & &  + \int_{\hat\Omega} \varepsilon c \varphi \,\mathrm{d}\omega_{\hat\Omega}\,.
    \label{variation_2_part_2} 
\end{eqnarray}
}
%%%%%%%%%%%%%%%%%%%%
%%%%%%%%%%%%%%%%%%%%
%%%%%%%%%%%%%%%%%%%%
\section{Coupled Discretized Model of DBV and Excess TBV Fractions}

We apply the Ritz--Galerkin method \cite{braess2007finite} for the discretization of $c$ and $\tilde q$ with trial function space defined by 
\begin{equation}
{\mathbb{V}_h} = \text{span}\{{\varphi}^1, \ldots, {\varphi}^m\} \subset {\mathbb{V}} 
\end{equation}
The discretization error is assumed to be orthogonal to the solution. For this purpose, $ { c} \in \mathbb{V}_h$, $ \tilde{ c} \in \mathbb{V}_h $ and  $\tilde q \in \mathbb{V}_h $ are approximated by linear combinations of the same set of nodal basis functions $\{\varphi^j\}_{j=1}^m$, which are piecewise linear and satisfy the Kronecker delta condition~\cite{cho-2007-kronecker}. The expressions for ${c}$ and $\tilde q$ at the nodal points of the finite element mesh in \(\hat{\Omega}\) are represented as:
\begin{align*}
{ c}({\bf x}; t) & = \sum_{i=1}^{m} c_i(t) {\varphi}^i({\bf x})\,, 
\\
\tilde{q}({\bf x}; t) & = \sum_{i=1}^{m} \tilde{q}_i(t) {\varphi}^i({\bf x})\,,
\end{align*}
where $ c_i $ and $ \tilde{q}_i $ are the nodal values of excess TBV and DBV fractions at each node $ { x}_i $ of the mesh, respectively, for $ i = 1, \ldots, m $. By utilizing the same basis functions for both $c$ and $\tilde{q}$, we maintain consistency in the numerical treatment of these related distributions within the discretized domain. Let's denote the coordinate vectors of the discretized variables $c$, $\bar c$, and $\tilde{q}$ by ${\bf c} = (c_1, c_2, \ldots, c_m)$, $\tilde{\bf c} = {\bf c} + \bar{\bf c}$ with $\bar{\bf c} = (\bar{c}_1, \bar{c}_2, \ldots, \bar{c}_m)$, and $\tilde{\bf q}$ $=$ $(\tilde{q}_1, \tilde{q}_2, \ldots, \tilde{q}_m)$, respectively. 
%%%%%%%%%%%%%%%
%%%%%%%%%%%%%%%
%%%%%%%%%%%%%%%
%\subsection{Discretized Volumetric Deoxygenated Hemoglobin  Concentration}
%\label{sec:discretized_blood_concentration}
By employing the implicit Euler method, the variational forms of \eqref{q_fick_variational_g} and \eqref{variation_2_part_2} can be respectively expressed as a coupled pair in the following form:
{\setlength\arraycolsep{2 pt} 
\begin{eqnarray}
\Big( {\bf M} + \Delta t \, ({\bf S} + {\bf T}  + {\bf G}) \Big)  \, \tilde{\bf q}^{(k+1)} 
 & = &  {\bf M} \tilde{\bf q}^{(k)} + \Delta t \,  {\bf b}^{(k)}  
 \label{model:DBV} 
 \\
  \Big(  {\bf M} + \Delta t \, ({\bf T} +   {\bf G})  \Big) {\bf c}^{(k+1)} & = & {\bf M} {\bf c}^{(k)} + \Delta t \,  {\bf s}\,,
  \label{model:TBV}
\end{eqnarray}
}
where the matrices and vectors are given by 
% {\setlength\arraycolsep{2 pt}\begin{eqnarray}
%    ({\bf M_{\tilde{\bf z}}})_{i, j} & = & \int_{\hat{\Omega}} \sum_{k = 1}^m ({\bf z}_k + {\bf z_k^{(0)}}) \, \varphi^k \varphi^i \varphi^j \, d \omega_{\hat{\Omega}} \\
%    ({\bf M_{{\bf y}}})_{i, j}  & = & \int_{\hat{\Omega}} \sum_{k = 1}^m {\bf y}_k \, \varphi^k \varphi^i \varphi^j \, d \omega_{\hat{\Omega}} \\
%     ({\bf K_{{\bf z}}})_{i, j}  & = & \int_{\hat{\Omega}} \delta \sum_{k = 1}^m {\bf z}_k \, \varphi^k {\bf g} ( \nabla \varphi^i, \nabla  \varphi^j) \, d \omega_{\hat{\Omega}} \\
%     ({\bf K_{{\bf y}}})_{i, j}  & = & \int_{\hat{\Omega}} \delta \sum_{k = 1}^m {\bf y}_k \, \varphi^k {\bf g} ( \nabla \varphi^i, \nabla  \varphi^j) \, d \omega_{\hat{\Omega}} 
%     \\
%     {\bf b}_i & = &  \int_{\hat\Omega} (\alpha \zeta + \upsilon h \sum_{k = 1}^m ({\bf z}_k + {\bf z_k^{(0)}}) \, \varphi^k) \varphi^i \,\mathrm{d}\omega_{\hat\Omega}
% \end{eqnarray}}
%%%%%%%%%%%%%%%%%%%
\begin{align*}
   M_{ij} & =  \int_{\hat{\Omega}} \,\varphi^i \varphi^j \, \mathrm{d} \omega_{\hat{\Omega}} 
   \\ 
T_{ij}  & =  \int_{\hat{\Omega}} \varepsilon  \,\varphi^i \,\varphi^j \, \mathrm{d} \omega_{\hat{\Omega}} 
    \\
    S_{ij}  & =  \int_{\hat{\Omega}} \upsilon  \,\varphi^i \,\varphi^j \, \mathrm{d} \omega_{\hat{\Omega}}  \\
    G_{ij}  & =  \int_{\hat{\Omega}} \delta   \,\varphi^i_{,h} \,\varphi^j_{,h} \, \mathrm{d} \omega_{\hat{\Omega}} 
    \\
    {b}^{(k)}_i & =   \int_{ \mathrm{B}}   (1-h) f \varphi^i \,\mathrm{d} \omega_{\mathrm{B}} +  \int_{\hat\Omega}  h\, \upsilon \tilde{c}^{(k)} \varphi^i \,\mathrm{d} \omega_{\hat{\Omega}}  \, \\
        {s}_i & =   \int_{\mathsf{B}}  f \varphi^i \, \mathrm{d} \omega_{\mathrm{B}} + \int_{\hat\Omega} \alpha_{, t}  \,\varphi^i \,\mathrm{d}\omega_{\hat\Omega} 
\label{vf_bc}
\end{align*}
\label{sec:discretized_blood_concentaration}
Due to the coupling between excess TBV and DBV, the k-th step approximation (\ref{model:TBV}) of excess TBV fraction has to be evaluated before the (k+1)-th step estimate (\ref{model:DBV}) of DBV fraction can be obtained.

%%%%%%%%%%%%%%%%%

% \subsubsection{Diagonal approximation}
% \label{sec:discretized_blood_flow}

% To simplify the evaluation of  (\ref{eq:g_matrix_form}), we approximate  matrices ${\bf M}_{\tilde{\bf c}}$, ${\bf M}_{{\bf q}}$, ${\bf G}_{\tilde{\bf c}}$, and ${\bf G}_{{\bf q}}$ via the products $ {\bf D}_{\tilde{\bf c}} {\bf M}$, $ {\bf D}_{{\bf q}} {\bf M}$, ${\bf D}_{\tilde{\bf c}} {\bf G}$, and ${\bf D}_{{\bf q}} {\bf G}$,  where ${\bf D}_{\tilde{\bf c}}$ and ${\bf D}_{{\bf q}}$ are diagonal matrices with diagonal entries given by 
% \begin{align*}
%    ( {\bf D}_{\tilde{\bf c}})_{ii} & =  \tilde{c}_i
%    \\ 
% ({\bf D}_{{\bf q}})_{ii}  & =   {q}_i
% \end{align*}
% for all $i=1,\ldots, m$. Substituting these approximations in (\ref{eq:g_matrix_form}) yields: 
% \begin{equation}
% \begin{aligned}
% {\bf D}_{\tilde{\bf c}} \Big( \big(1 + \Delta t \,  (\upsilon + \varepsilon) \big)  {\bf M} &+ \Delta t \,  {\bf G} \Big)  \, {\bf q}^{(k+1)} 
%  = {\bf D}_{\tilde{\bf c}} {\bf M} \, {\bf q}^{(k)}
% \\
% & + {\bf D}_{\bf q} \, {\bf M} \, \Big({\tilde{\bf c}}^{(k)} - {\tilde{\bf c}}^{(k+1)} \Big)
%  + \Delta t \,  \Big({\bf b}   - {\bf D}_{\bf q} \, {\bf G} \, {\tilde{\bf c}}^{(k)} \Big) \,.
% \end{aligned}
% \label{eq:g_matrix_form_2}
% \end{equation}

%%%%%%%%%%%%%%%

\section{Head Model and Vessel Segmentation via 7T MRI data}

To test the performance of the coupled model presented by equations (\ref{model:DBV}) and (\ref{model:TBV}), we utilized a multi-compartment head model and vessel segmentation \cite{samavaki2023pressure}, derived from high-resolution imaging data and advanced image processing techniques. The head model was developed using the open sub-millimeter precision 7 Tesla MRI dataset of CEREBRUM-7T \cite{svanera2021cerebrum}. 

The use of 7T MRI data is crucial for accurate vessel reconstruction due to its high spatial resolution, which allows the identification of vessels with diameters larger than the voxel size \cite{fiederer2016}. Lower field MRI systems (e.g., 3T) may not capture small cerebral vessels, leading to incomplete vascular models. In particular, high-resolution data is essential for detailed computational models that simulate blood flow, where precise vessel geometry is necessary to model the dynamic effects of cerebral blood flow.

The present 7T MRI dataset was acquired using a Siemens 7T Terra Magnetom MRI scanner with a 32-channel head coil; see  \cite{svanera2021cerebrum, our-7T-dataset}. The images were obtained using the MP2RAGE sequence at an isotropic resolution of 0.63 mm$^3$. Data collection took place at the Imaging Centre of Excellence at the Queen Elizabeth University Hospital, Glasgow, UK. A total of 142 brain volumes were collected as reconstructed DICOM images.

%%%%%%%%%%%%%%%%%%%%%%%%
%%%%%%%%%%%%%%%%%%%%%%%%
%%%%%%%%%%%%%%%%%%%%%%%%
\subsection{Reconstruction of Vessels Using Frangi's Vesselness Algorithm}

\begin{figure}[h!]
\centering
\begin{footnotesize}
\begin{minipage}{4.25cm}
\centering
\includegraphics[width=4cm]{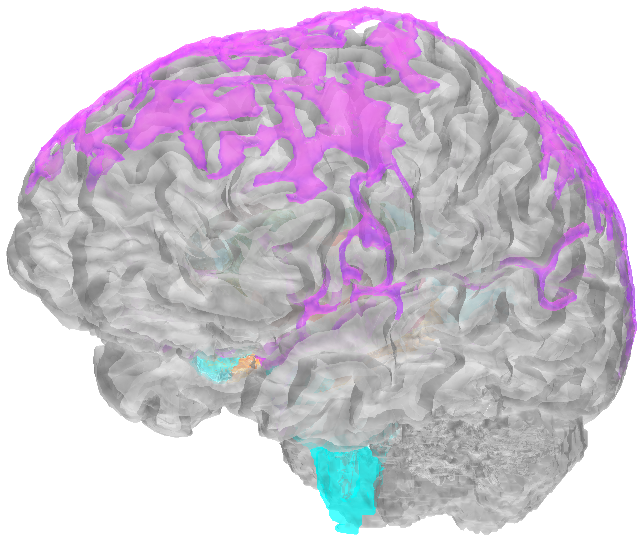} \\
\end{minipage}
\begin{minipage}{4.25cm}
\centering
\includegraphics[width=4cm]{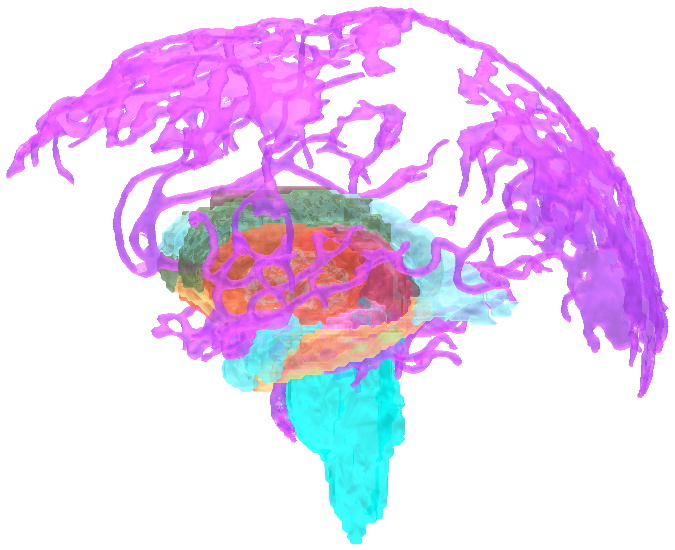} \\
\end{minipage} \\ \vskip0.1cm
\includegraphics[width=4cm]{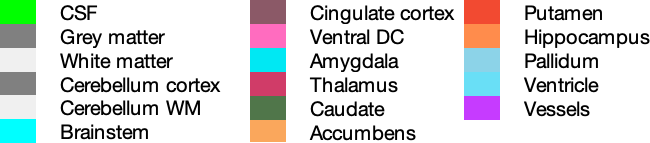}
\end{footnotesize}
\caption{The vessel segmentation of this study together with the grey matter (GM) layer (left) and subcortical compartments (right) of the head model. \label{fig:head_model_and_vessel_segmentation}}
\end{figure}

In this study, the vessel segmentation was carried out using Frangi's Vesselness algorithm \cite{frangi1998multiscale,choi2020cerebral}, which is specifically designed to enhance tubular structures in medical images. Initially, Frangi's algorithm was applied to the preprocessed INV2 and T1w MRI slices separately, utilizing the implementation available in the Scikit-Image package \cite{van2014scikit}. The results were then aggregated by superposing the filtered slices and binarizing the mask with a user-defined threshold. Noise reduction was achieved by aggregating results from sagittal, axial, and coronal slices, ensuring that voxels detected as vessels in two or more views were retained in the final vessel mask.  The vessel segmentation is illustrated in Figure \ref{fig:head_model_and_vessel_segmentation}.

%%%%%%%%%%%%%%%%%%%%%%
%%%%%%%%%%%%%%%%%%%%%
%%%%%%%%%%%%%%%%%%%%%
\subsection{Segmentation of Head Compartments}

The segmentation of the other brain model compartments was executed using the FreeSurfer software suite \cite{fischl2012freesurfer}, in combination with the FieldTrip interface \cite{oostenveld2011} for the SPM12 surface extractor \cite{ashburner2014spm12}. The segmentation included 17 brain compartments: cerebrospinal fluid (CSF), grey matter, white matter, cerebellum cortex, cerebellum white matter, brainstem, cingulate cortex, ventral diencephalon, amygdala, thalamus, caudate, accumbens, putamen, hippocampus, pallidum, ventricles, and vessels. The segmentation methods are summarized in Table \ref{tab:segmentation}. The microvessel length densities $\xi$ in Table \ref{tab:microvessel_length_densities} in different compartments of the segmentation were estimated according to the median values found in \cite{kubikova2018numerical}.

\begin{table}[h!]
    \caption{Segmentation methods applied for the blood vessel, grey matter (GM), white matter (WM), and cerebrospinal fluid (CSF) compartments constituting the brain model of this study. The software applied included Frangi's Vesselness algorithm \cite{frangi1998multiscale}, FreeSurfer software suite \cite{fischl2012freesurfer}, and FieldTrip's SPM12 head and brain segmentation pipeline.
    \label{tab:segmentation}}
    \begin{footnotesize}
    \begin{center}
    \begin{tabular}{ll}
        \toprule
        Compartment & Method  \\
        \hline
        Blood vessels & Vesselness   \\
        Cerebral, cerebellar \& subcortical\textsuperscript{$\ast$} GM & FreeSurfer  \\
        Cerebral, cerebellar \& subcortical\textsuperscript{$\dag$} WM & FreeSurfer \\
        Ventricles (CSF) & FreeSurfer  \\
        CSF & FieldTrip/SPM12 \\
        \bottomrule
    \end{tabular} 
    \end{center} 
$\ast$ Brainstem, Ventral diencephalon, Amygdala, Thalamus, Caudate, Accumbens, Putamen, Hippocampus, Pallidum.\\
$\dag$ Cingulate cortex.
    \end{footnotesize}
\end{table}

\begin{table}[ht!]
\begin{footnotesize}
\begin{center}
\caption{Microvessel length densities across different grey and white matter compartments according to the median values found in~\cite{kubikova2018numerical}. Some compartments have identical density values, reflecting their structural similarities or broader classification into regions with similar vascular traits. }
\label{tab:microvessel_length_densities}
\mbox{} \\
\begin{tabular}{lc}
\toprule
{Brain Compartment} & {Length Density $\xi$ ($\text{m}^{-2}$)} \\
\hline
Cerebral GM       & $2.4 \times 10^8$ \\
Cerebral WM       & $1.4 \times 10^8$ \\
Cerebellar GM      & $3.0 \times 10^8$ \\
Cerebellar GM     & $1.0 \times 10^8$ \\
Subcortical\textsuperscript{$\ast$} GM    & $3.3 \times 10^8$ \\
Subcortical\textsuperscript{$\dag$}  WM    & $1.5 \times 10^8$ \\
Brainstem                   & $2.9 \times 10^8$ \\
\bottomrule
\end{tabular}
\end{center}
$\ast$ Ventral diencephalon, Amygdala, Thalamus, Caudate, Accumbens, Putamen, Hippocampus, Pallidum.
\\
$\dag$ Cingulate cortex.
\end{footnotesize}
\end{table}

%%%%%%%%%%%%%%%%%%%%%%
%%%%%%%%%%%%%%%%%%%%%%
%%%%%%%%%%%%%%%%%%%%%%
\subsection{Numerical experiments}

For our numerical experiments, we performed three simulation runs utilizing the coupled model given by (\ref{model:DBV}) and (\ref{model:TBV}) to analyze DBV and excess TBV fractions, respectively. The time interval covered by the models was 21 seconds, with a time step length of 0.25 seconds used in the iterative computation runs. Other central model parameters can be found in Table \ref{tab:physical_parameters}. An estimate for the excess TBV fraction was first found by solving (\ref{model:TBV}) with a boundary flux condition following a static approximation obtained via the pressure-Poisson equation (PPE) as described in \cite{samavaki2023pressure} and the parameters therewithin. 

\begin{figure}[h!]
\centering
\begin{footnotesize}
\begin{minipage}{4.25cm}
\centering
\includegraphics[height=3cm]{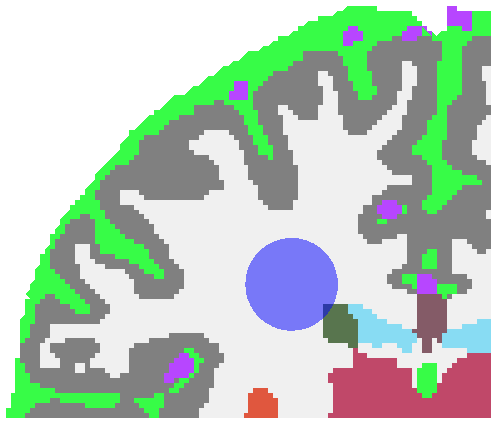} \\
WM-ROI \\ Coronal
\end{minipage}
\begin{minipage}{4.25cm}
\centering
\includegraphics[height=3cm]{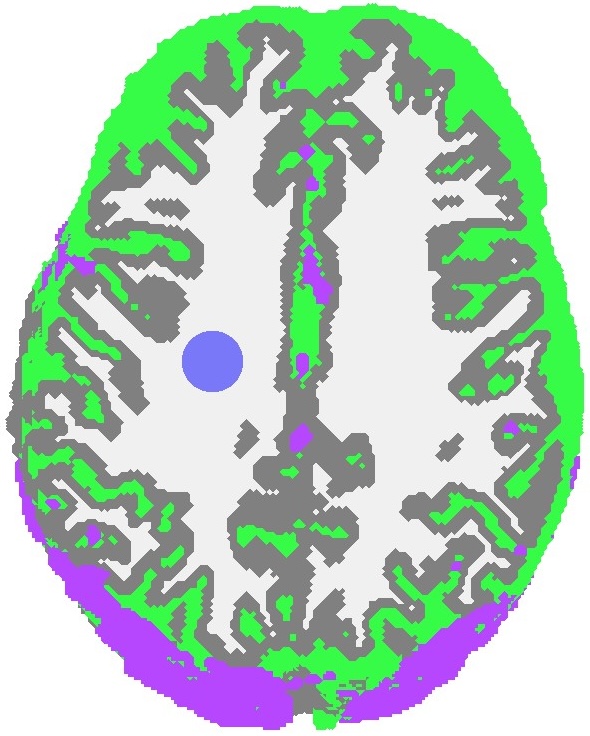} \\
WM-ROI \\ Axial
\end{minipage} \\ \vskip0.2cm
\begin{minipage}{4.25cm}
\centering
\includegraphics[height=3cm]{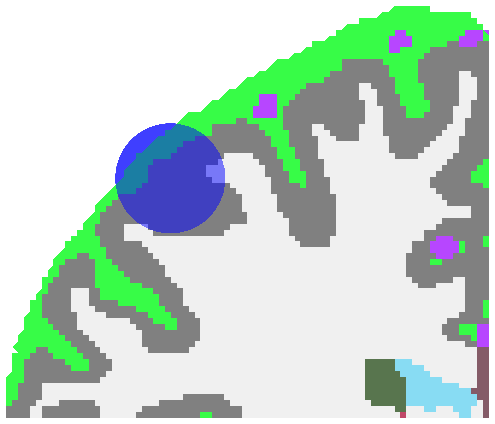} \\ 
GM-ROI \\
Coronal
\end{minipage}
\begin{minipage}{4.25cm}
\centering
\includegraphics[height=3cm]{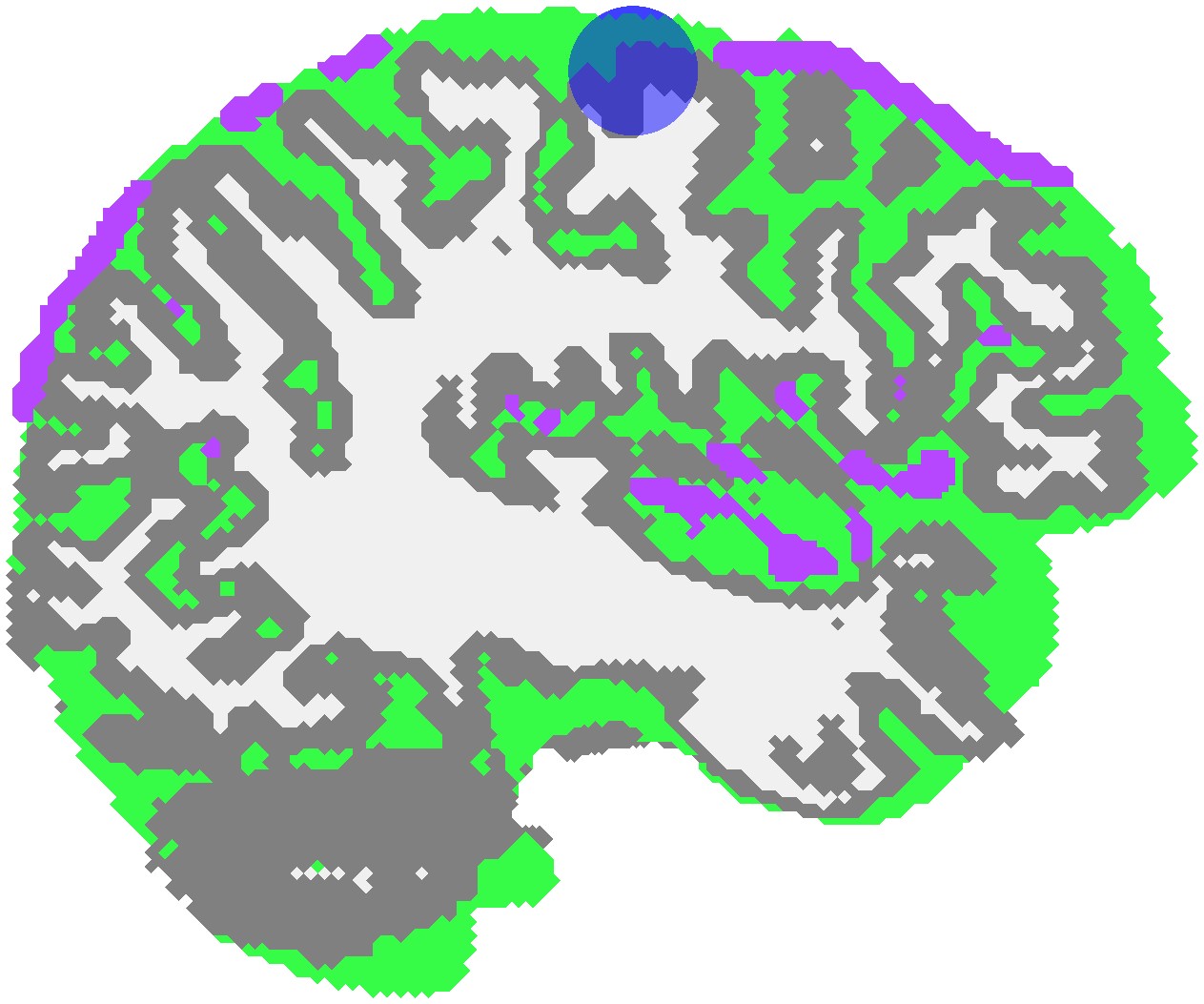} \\ 
GM-ROI \\
Sagittal
\end{minipage} \\ \vskip0.2cm
\begin{minipage}{4.25cm}
\centering
\includegraphics[height=3cm]{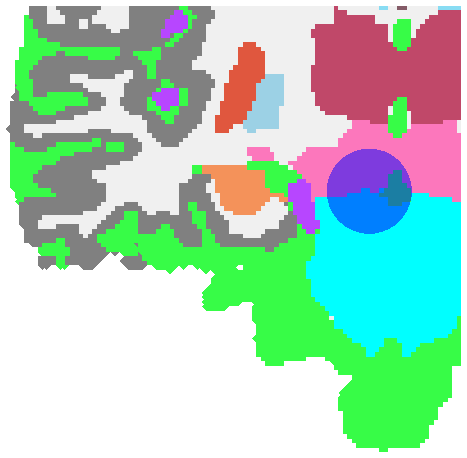} \\
BS-ROI \\ Coronal
\end{minipage}
\begin{minipage}{4.25cm}
\centering
\includegraphics[height=3cm]{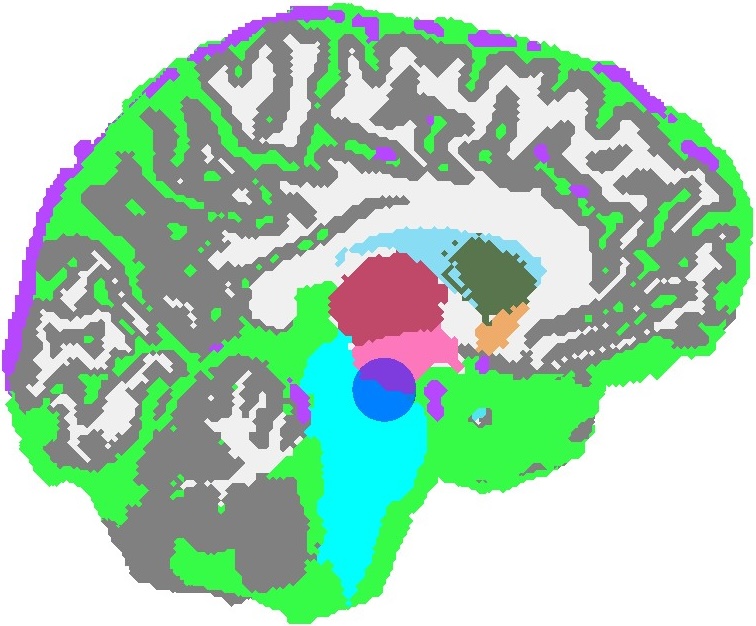} \\
BS-ROI \\ Sagittal
\end{minipage} \\ \vskip0.2cm
\begin{minipage}{6cm}
\centering
\includegraphics[width=4cm]{legend2.png}
\end{minipage}
\end{footnotesize}
\caption{Three different regions of interest (ROIs), each with 14-mm-diameter, were placed in white matter (WM), grey matter (GM), and the brainstem (BS). The ROIs share center points within the same coronal slice used for visualization on the left.  
\label{fig:spatial_rois} }
\end{figure}

We analyzed the spatial and temporal distributions of DBV and excess TBV fractions within three distinct regions of interest (ROIs), each with a 14-mm diameter. These ROIs are located in the cerebral white matter (WM-ROI), grey matter (GM-ROI), and the subcortical region of the brainstem (BS-ROI), which partially overlaps with the ventral diencephalon compartment. The locations of these regions are illustrated in Figure \ref{fig:spatial_rois}. These ROIs have different densities of small blood vessels and varying connections to arteries which  affects both the background $\bar{ c}$ and excess ${ c}$ blood volume fraction distribution in these areas. This variation allows us to study how local changes in $\bar{ c}$ and ${ c}$ influence our model. 

As a software platform, we used the Matlab-based Zeffiro Interface (ZI) \cite{he2019zeffiro}. The volume of the head segmentation was discretized into a tetrahedral FE mesh with an approximate resolution of 1 mm. Specifically, the region $\Omega$ consisted of 0.15 million nodes and 0.54 million tetrahedra, whereas the region $\hat{\Omega}$ contained 2.4 million nodes and 11 million tetrahedra. For generating the FE mesh, we applied a Dell Precision 5820 Tower Workstation equipped with 256 GB of random-access memory (RAM), an Intel Core i9-10900X 3.70 GHz (10 cores and 20 threads) central processing unit (CPU), and two Quadro RTX 8000 GPUs with 48 GB of GPU RAM in each. %The iterations of the coupled model of (\ref{model:DBV}) and (\ref{model:TBV}) were performed using a laptop with 16 GB of RAM and 8-core Apple M1 CPU without GPU acceleration.

The results in this work can be fully reproduced using ZI's \cite{pursiainen_zeffiro_august_2024}, version 6.0.2. A head model can be  generated according to ZI's standard head model generator, which first imports a surface segmentation and then generates an FE mesh based on that \cite{galaz2023multi}. Instructions for generating such a head model can be found in ZI's repository. All the solvers applied in this study are included in ZI's NSE tool plugin. The visualizations of this study have been created using ZI's visualization tool.

\begin{table}[h!]
\centering
\caption{Parameters of numerical Simulations: The neural drive parameter $\zeta$ has been adjusted so that the current time-mollified balloon model matches a closely matches $20 \%$ increase in local blood volume, as observed in experimental data from the cat cortex \cite{kim2011temporal}. The signal decay rate parameter $\kappa$ and the flow-dependent elimination constant parameter $\gamma$ were chosen based on values from \cite{friston2003dynamic}. The fraction oxygen consumption rate $\upsilon = \eta \rho^{-1} \psi^{-1}$ has been obtained from the molecular rate 0.4 $s^{-1}$ used in \cite{buxton2004modeling}. Molality $\psi$ of hemoglobin in blood is as in \cite{aquino2014spatiotemporal}. Values for blood density $\rho$,  viscosity $\mu$, total cerebral blood flow (CBF) $Q$, relative pressure drop in arterioles $\vartheta$, and the diameters of arterioles, capillaries, and venules ($D_a$, $D_c$, and $D_v$, subtracting the total wall thicknesses of 2.0E-5, 2.0E-06, and 2E-6, respectively), along with their corresponding total area fractions $\gamma_a$, $\gamma_c$, and $\gamma_v$, can be found in textbooks \cite{caro2012mechanics, tu2015human}. The reference level for the oxygenated blood volume (OBV) was selected under the assumption that it is in the capillary circulation, corresponding to approximately half of the total oxygen extraction fraction found in white matter  \cite{ito2023oxygen}. The reference pressures were selected based on \cite{gabriel1996compilation} and \cite{blanco2017blood}, respectively. Gravitational acceleration was adjusted to its average value and aligned with the z-axis.}
\label{tab:physical_parameters}
\resizebox{0.95\linewidth}{!}{
    \begin{tabular}{llr}
    \toprule
 Property & Parameter    & Value\\
        \midrule
Time interval length & $T$ & $21 \ \text{s}$
\\
Time step length & $\Delta t$ &  $0.25 \text{s}$ 
%\\
%{\color{red}  Arteriole length } & {\color{red}$L$} & {\color{red}$0.4 \text{mm}$}
\\
%{\color{red}  Coupling constant coefficient  } & {\color{red}$\tau$} & {\color{red}$ ?$}
%\\
 Average diastolic pressure  & $\overline{p}$ & $75 \text{mmHg}$   
\\
Neural drive & $\zeta$ & $1 \text{s}^{-2}$ 
\\
Signal decay rate & $\kappa$ & $0.65 \text{s}^{-1}$ 
\\
Flow-depended elimination constant & $\gamma$ & $0.41 \,  \text{s}^{-2}$  
\\
%Volumetric oxygen consumption rate & $\upsilon$ & $0.023 \text{s}^{-1}$
%\\ 
Molecular oxygen consumption rate & $\eta$ & $0.4 \, \text{s}^{-1}$
\\ 
Molality of hemoglobin in blood & $\psi$ & $0.0018 \, \text{mol kg}^{-1}$ 
\\
Blood density                       & $\rho$          & $1050 \,  \text{kg m}^{-3}$ \\
Blood viscosity                     & $\mu$      & $4.0 \,  \times 10^{-3} \text{Pa s}$ 
\\
Total cerebral blood flow           & $Q$       & $750 \, \text{ml min}^{-1}$     \\
Arteriole diameter                  & $D_a$           & $1.0   \times 10^{-5} \, \text{m}$  
\\
Capillary diameter      & $D_c$   & $7.0 \times 10^{-6} \, \text{m}$  
\\
Venule diameter                     & $D_v$           & $1.8 \times 10^{-5} \, \text{m}$   
\\
%{\color{red} Microvessels in cerebral GM} & {\color{red} $\xi$ }& {\color{red} $2.4\times 10^{-8} \,  \text{m}^{-2}$ } 
%\\
%{\color{red} Microvessels in cerebral WM} & {\color{red} $\xi$} & {\color{red}$ 1.4\times 10^{-8} \,  \text{m}^{-2}$ }
%\\
%{\color{red} Microvessels in cerebellar GM} & {\color{red} $\xi$} & {\color{red} $3.0\times 10^{-8} \,   \text{m}^{-2}$ } 
%\\
%{\color{red} Microvessels in cerebellar WM} & {\color{red} $\xi$} & {\color{red} $1.0\times 10^{-8} \,  \text{m}^{-2}$ }  
%\\
%{\color{red} Microvessels in subcortical WM }& {\color{red} $\xi$} & {\color{red} $1.5\times 10^{-8} \,  \text{m}^{-2}$ }  
%\\
%{\color{red} Microfibers in brainstem }&{\color{red}  $\xi$} & {\color{red} $2.9\times 10^{-8} \,  \text{m}^{-2}$ }  
%\\
Arteriole total area fraction & $\gamma_a$ &  30 \%
\\
Capillary area fraction & $\gamma_c$ &   40 \%
\\
Venule area fraction & $\gamma_v$ & 30 \%
%\\
%{\color{red} Resting state blood volume fraction } & {\color{red}$\overline{c}$} & {\color{red}$? \text{?}$}
%\\
%{\color{red} Reference deoxygenated blood volume fraction } & {\color{red}$\overline{q}$} & {\color{red}$? \text{?}$}
\\
Reference oxygenated blood volume fraction & $h$ &  $85\%$ \\
Relative pressure decay in arterioles        & $\vartheta$        & $70\%$            
\\
Microvessel cross-sectional area expansion               & $\theta$      & $0.21$                \\ 
 Gravitation (z-component) & $f^{(3)}$ &  \, -9.81\,\text{m s\textsuperscript{-2}} \\
 \bottomrule
    \end{tabular}
} % end resizebox
\label{tab:parameters}
\end{table}

%%%%%%%%%%%%%%

%%%%%%%%%%%%%%%%%%
%%%%%%%%%%%%%%%%%%
%%%%%%%%%%%%%%%%%%
\section{Results}

\begin{figure}[h!]
\centering
\includegraphics[height=2.5cm]{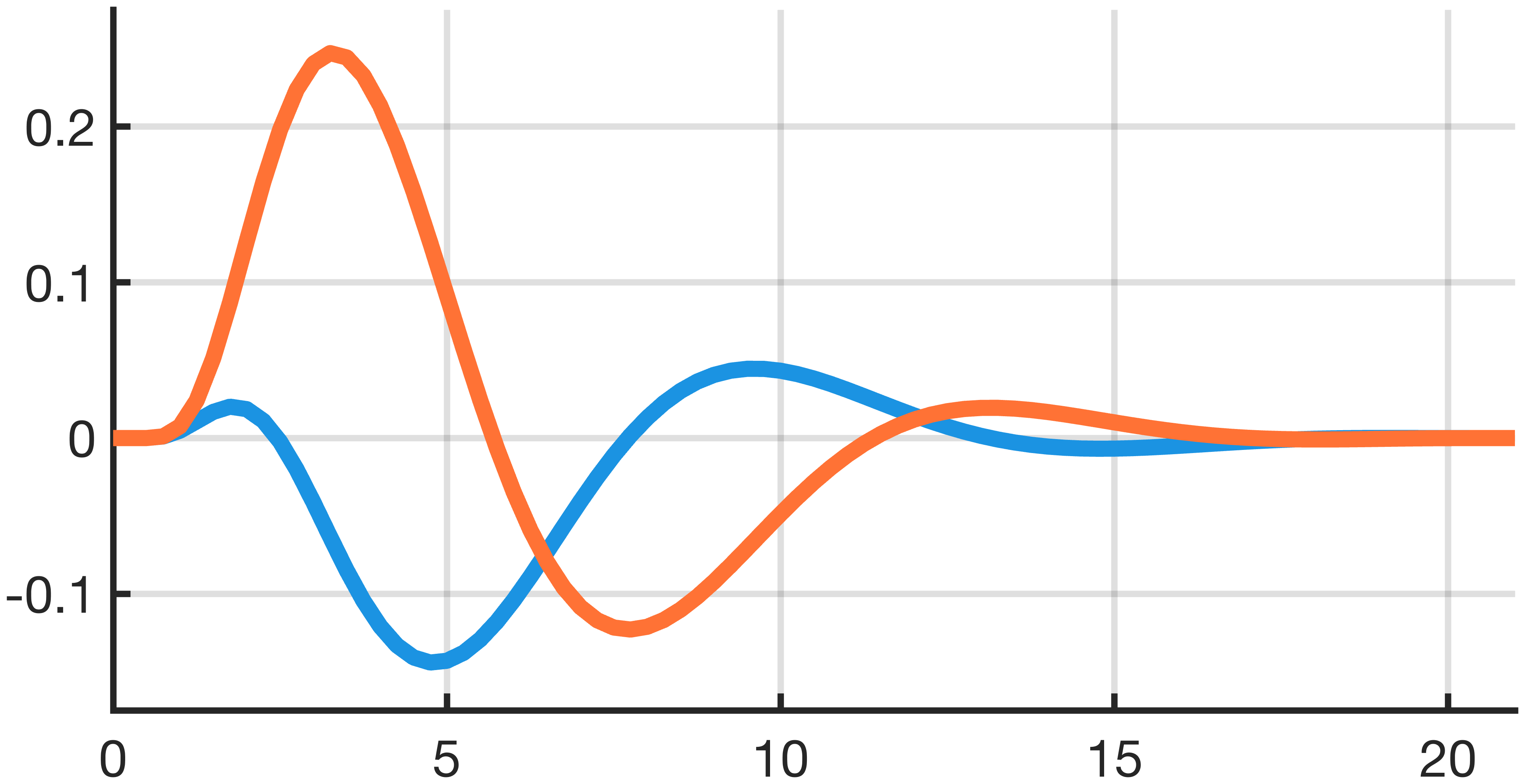}
\caption{The mollified approximation of hemodynamic response function $\alpha(t)$ (red) and its time derivative $\alpha_{,t} (t)$ (blue) are shown as a result of the time-mollified balloon model, i.e., a mollified solution of (\ref{DHOM}) as explained in Section \ref{sec:balloon_model}, with parameters selected as described in Table \ref{tab:physical_parameters}. An amplitude of approximately 20 \% has been measured, e.g., in the cat cortex \cite{kim2011temporal}.  
The horizontal axis represents time $t$ in seconds, while the vertical axis represents the value of  $\alpha(t)$.\label{fig:balloon_model}}
\end{figure}

\begin{figure*}[h!]
\centering
\begin{footnotesize}
\begin{minipage}{5.5cm}
\centering
\includegraphics[height=2.5cm]{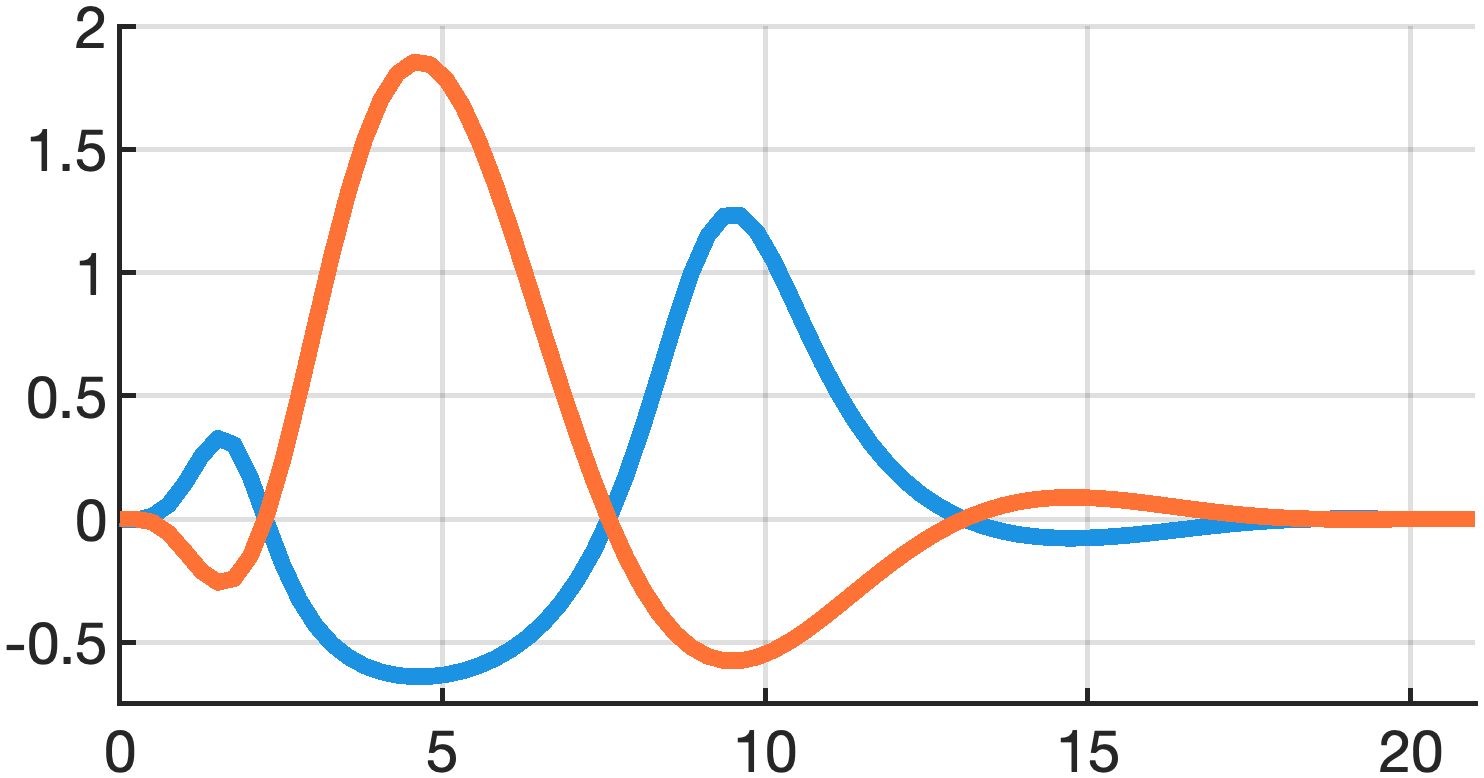} \\
WM-ROI \\
\end{minipage}
\begin{minipage}{5.5cm}
\centering
\includegraphics[height=2.5cm]{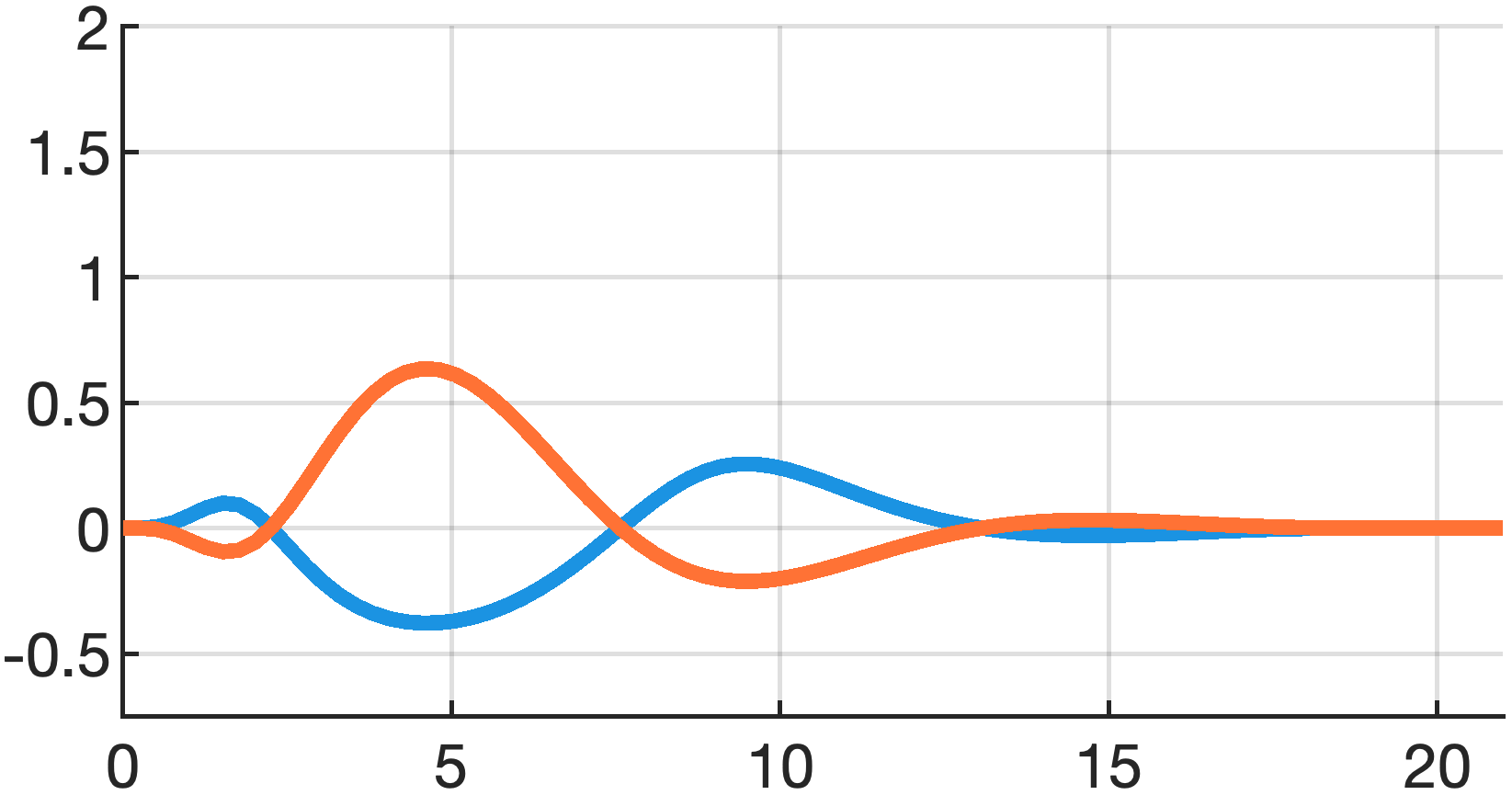} \\
GM-ROI \\
\end{minipage}
\begin{minipage}{5.5cm}
\centering
\includegraphics[height=2.5cm]{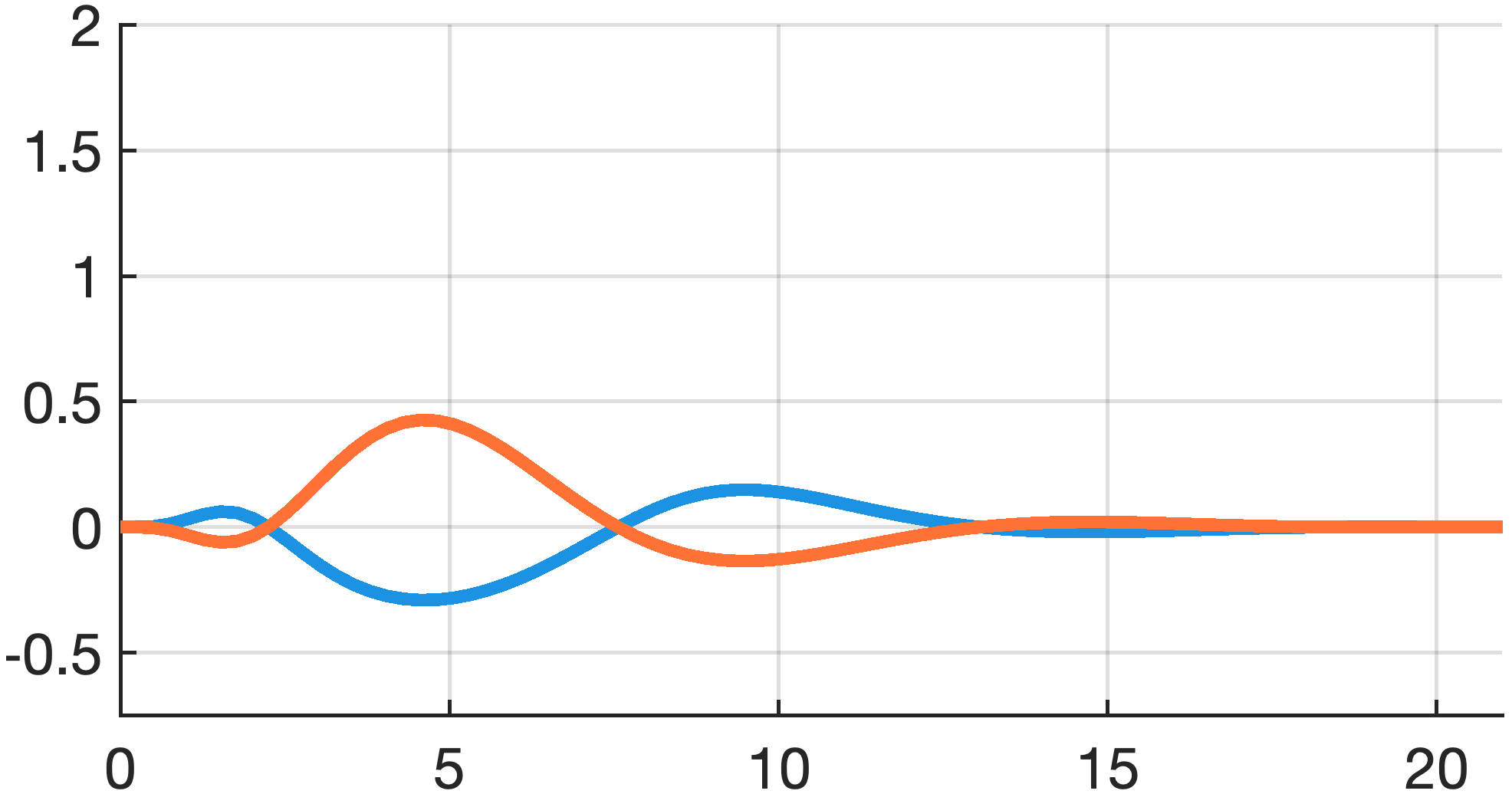} \\
BS-ROI \\
\end{minipage} \\
\end{footnotesize}
\caption{Time-evolution of the local mean (the spatial integral mean over the corresponding ROI for each time point)  relative oxygenated (red) and deoxygenated (blue) blood volume (OBV and DBV) fraction with respect to  three different regions of interest (ROIs) of this study. The  horizontal axis shows the time $t$ (s).  In each case, the neural activity is modelled as a point source in the center of the ROI.\label{fig:oxy_vs_deoxy_temporal}}
\end{figure*}

\begin{figure}[h!]
\centering
\begin{footnotesize}
\begin{minipage}{4.25cm}
\centering
\includegraphics[height=3cm]{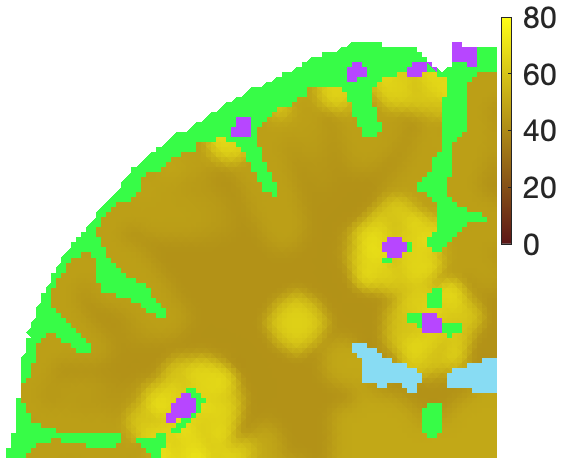} \\
WM-ROI \\ Coronal
\end{minipage}
\begin{minipage}{4.25cm}
\centering
\includegraphics[height=3cm]{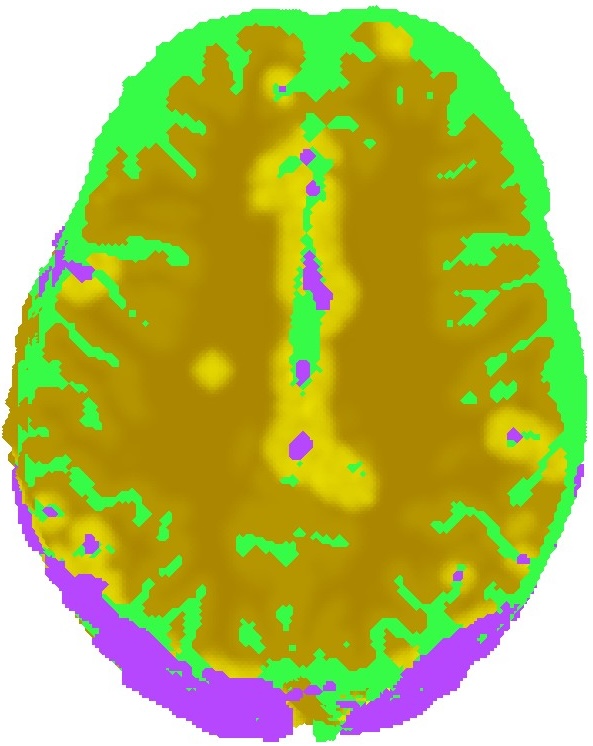} \\
WM-ROI \\ Axial
\end{minipage} \\ \vskip0.1cm
\begin{minipage}{4.25cm}
\centering
\includegraphics[height=3cm]{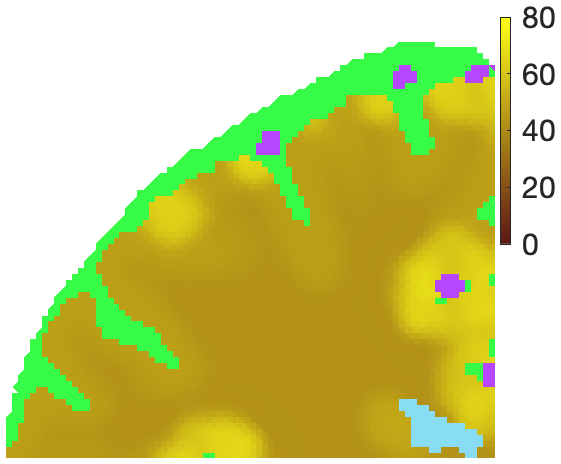} \\
GM-ROI \\ Coronal
\end{minipage}
\begin{minipage}{4.25cm}
\centering
\includegraphics[height=3cm]{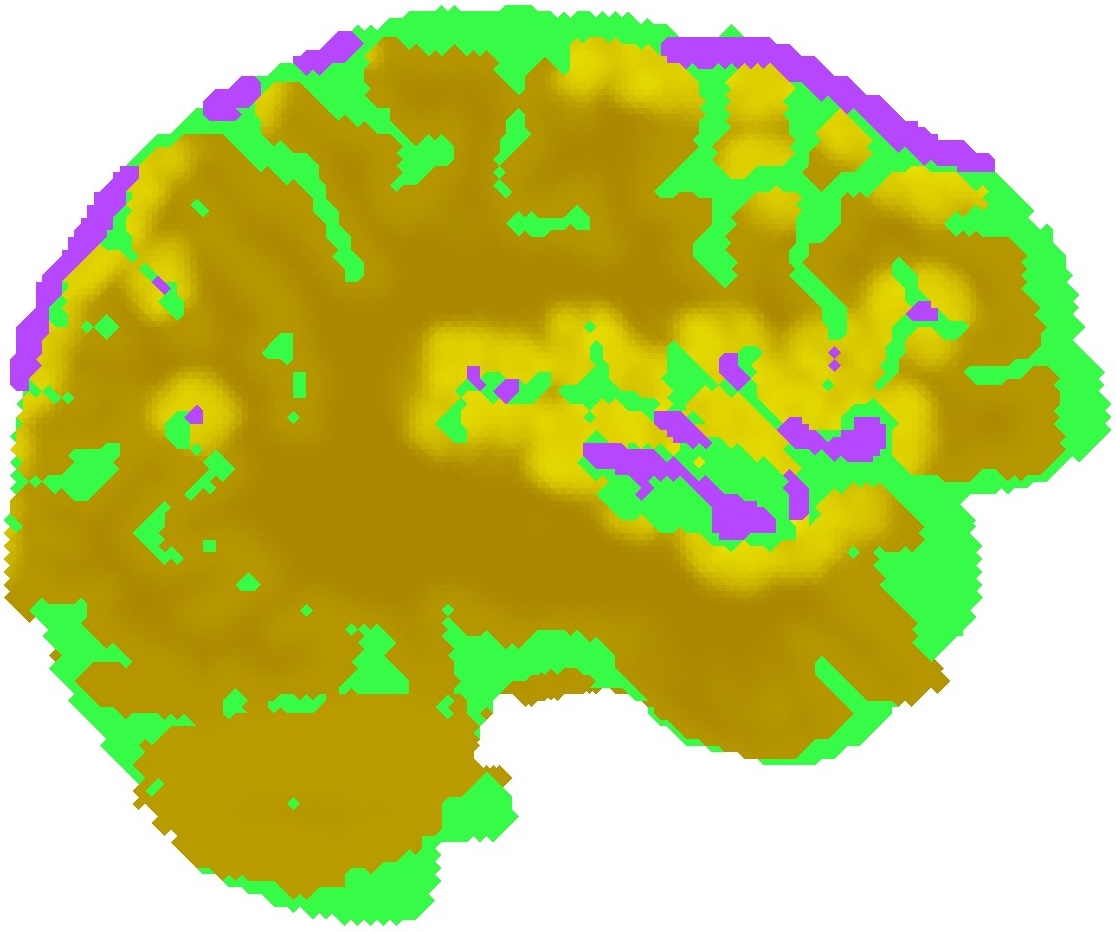} \\
GM-ROI \\ Sagittal
\end{minipage} \\ \vskip0.1cm
\begin{minipage}{4.25cm}
\centering
\includegraphics[height=3cm]{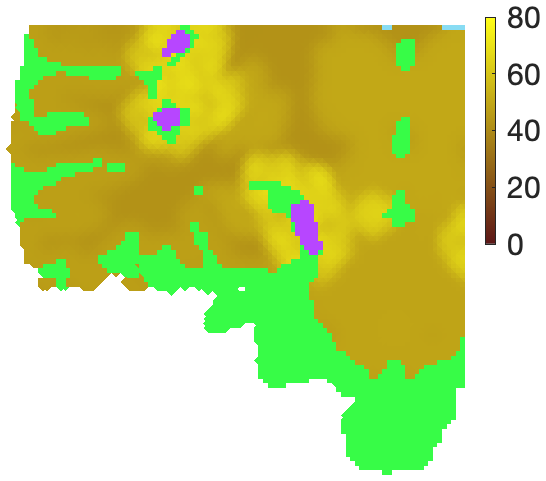} \\
BS-ROI \\ Coronal
\end{minipage}
\begin{minipage}{4.25cm}
\centering
\includegraphics[height=3cm]{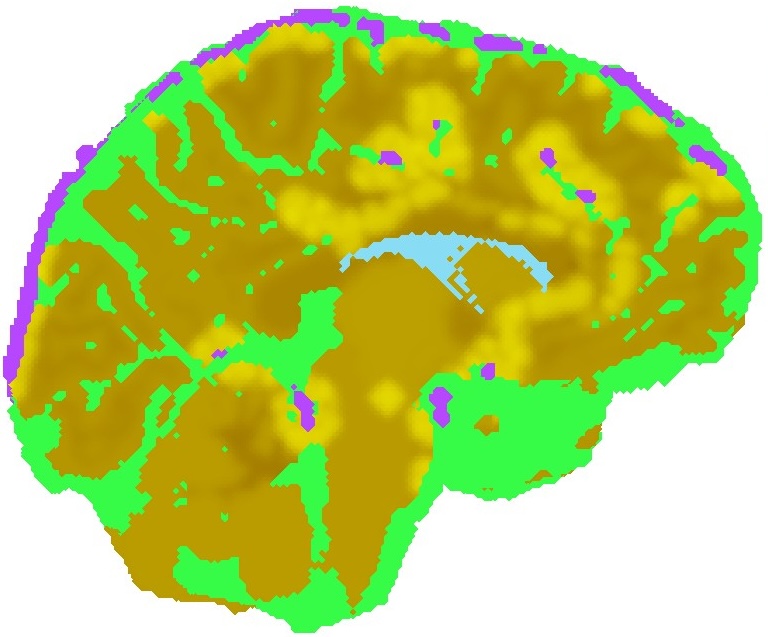} \\
BS-ROI \\ Sagittal
\end{minipage} 
\vskip0.1cm
\begin{minipage}{6cm}
\centering
\includegraphics[width=4cm]{legend2.png}
\end{minipage}
\end{footnotesize}
\caption{Total volumetric blood (TBV) fraction at $t = 5$ s on a decibel scale in three different cases, in which a point-like neural activity is placed in the center of WM-ROI, GM-ROI, and BS-ROI. The hemodynamic response can be recognized as it does not have a vessel in the center as a source of the excess concentration. The coronal slice of the visualization is shared by the ROIs' center points.  The color scale is logarithmic between 0 and 80 dB, where the upper limit corresponds to the maximum of the distribution. \label{fig:VB_spatial} }
\end{figure}

\begin{figure}[h!]
\centering
\begin{footnotesize}
\begin{minipage}{4.25cm}
\centering
\includegraphics[height=3cm]{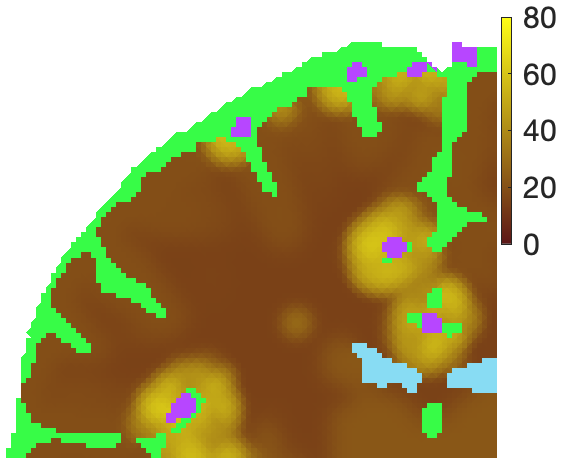} \\
WM-ROI \\ Coronal 
\end{minipage}
\begin{minipage}{4.25cm}
\centering
\includegraphics[height=3cm]{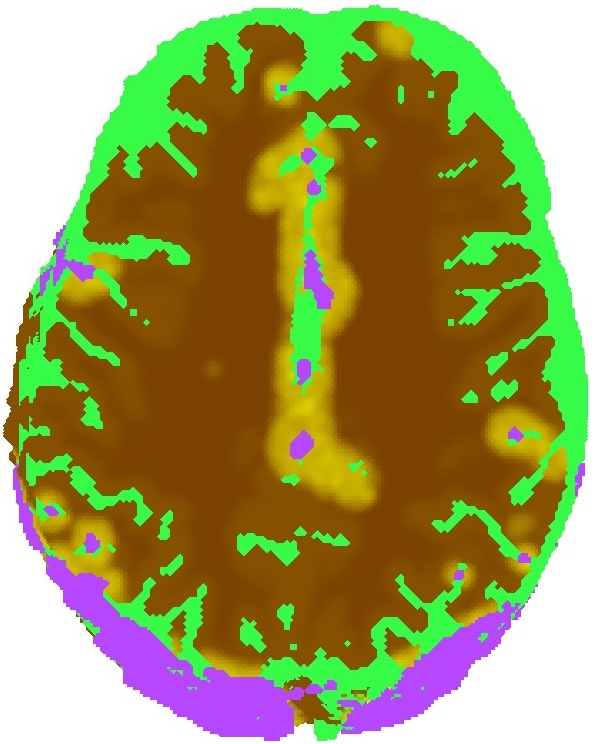} \\
WM-ROI \\ Axial 
\end{minipage}
\vskip0.1cm
\begin{minipage}{4.25cm}
\centering
\includegraphics[height=3cm]{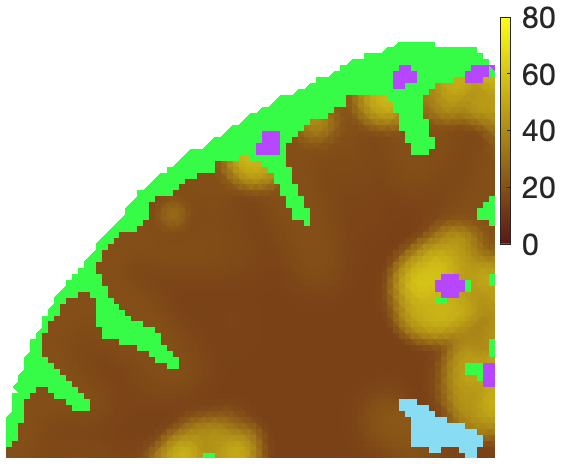} \\
GM-ROI \\ Coronal
\end{minipage}
\begin{minipage}{4.25cm}
\centering
\includegraphics[height=3cm]{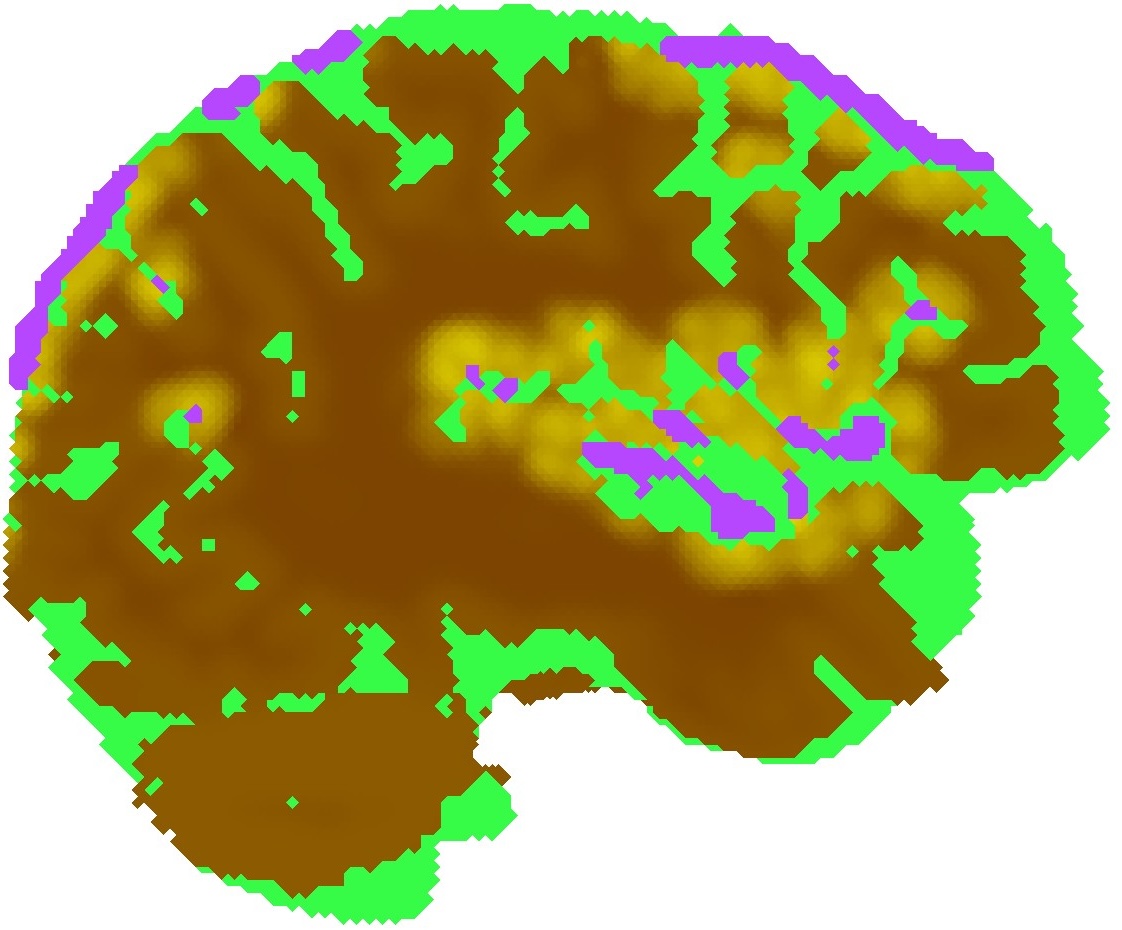} \\
GM-ROI \\ Sagittal
\end{minipage}
\vskip0.1cm
\begin{minipage}{4.25cm}
\centering
\includegraphics[height=3cm]{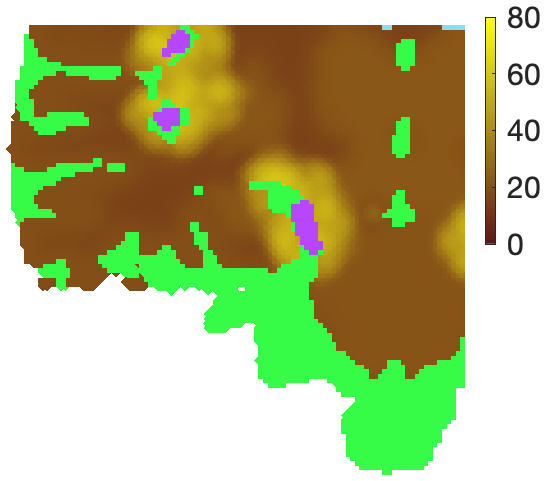} \\
BS-ROI \\ Coronal
\end{minipage}
\begin{minipage}{4.25cm}
\centering
\includegraphics[height=3cm]{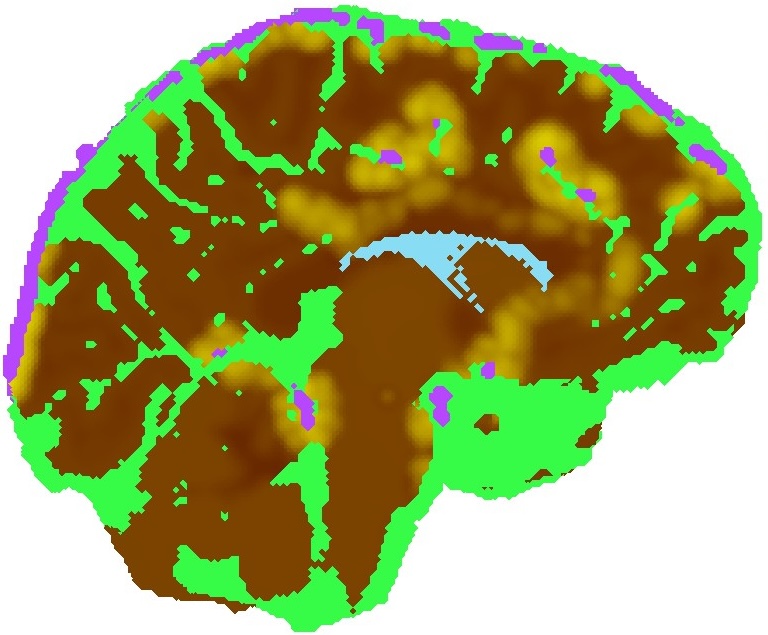} \\
BS-ROI \\ Sagittal
\end{minipage} \\
\vskip0.1cm
\begin{minipage}{6cm}
\centering
\includegraphics[width=4cm]{legend2.png}
\end{minipage}
\end{footnotesize}
\caption{Deoxygenated blood volume (DBV) fraction at $t = 5$ s on a decibel scale in three different cases, in which a point-like neural activity is placed in the center of WM-ROI, GM-ROI, and BS-ROI. The hemodynamic response can be recognized as it does not have a vessel in the center as a source of the excess concentration. The coronal slice of the visualization is shared by the ROIs' center points.  The color scale is logarithmic between 0 and 80 dB, where the upper limit corresponds to the maximum of the distribution. \label{fig:dHb_spatial} }
\end{figure}

Figure \ref{fig:balloon_model} shows the time-smoothed hemodynamic response function and its derivative, both derived from the balloon model (Section \ref{sec:balloon_model}). The resulting relative temporal changes in both OBV and DBV fractions within the three ROIs -- white matter (WM), grey matter (GM), and brainstem (BS) --  are illustrated in Figure \ref{fig:oxy_vs_deoxy_temporal}. These results indicate that local variations in DBV and TBV fractions significantly influenced the spatiotemporal characteristics of the hemodynamic response. The results demonstrate a localized increase in both DBV and TBV fractions around the site of neural activity, with a gradual spread to surrounding areas (Figures \ref{fig:VB_spatial} and \ref{fig:dHb_spatial}). In the WM-ROI, the highest TBV fraction increase is observed at $t=5$ s, with an approximate background level of 0.01 near the source of neural activity. The GM-ROI and BS-ROI, however, exhibit background levels of approximately 0.02 (as shown in Table \ref{tab:bloodflow}), suggesting that these regions are less responsive to localized neural activity than the WM-ROI. Our simulations specifically show, in the beginning of the observed time interval, that relative level of OBV fraction increases rapidly following local hemodynamic changes, and then diffuses outward according to Fick's law, highlighting the transport dynamics within the microvasculature (Figure \ref{fig:oxy_vs_deoxy_temporal}).  

%Spatially, table \ref{tab:bloodflow} presents information on TBV and DBV, their ratio $\tilde{q}/\tilde{c}$, and the magnitude of the perturbations in these distributions caused by neural activity within the ROIs for $t = 0.25$ s and $t = 5$ s.  This is reflected in the perturbation diameters reported at $t=5$ s, where the WM-ROI shows the largest TBV perturbation diameter of 8.49 mm, compared to 8.42 mm for GM-ROI and 7.96 mm for BS-ROI. Similarly, the DBV perturbation diameters are smaller, ranging from 4.05 mm in WM-ROI to 3.35 mm in BS-ROI.}

\begin{figure*}[h!]
\centering
\begin{footnotesize}
\begin{minipage}{5.5cm}
\centering
\includegraphics[height=2.5cm]{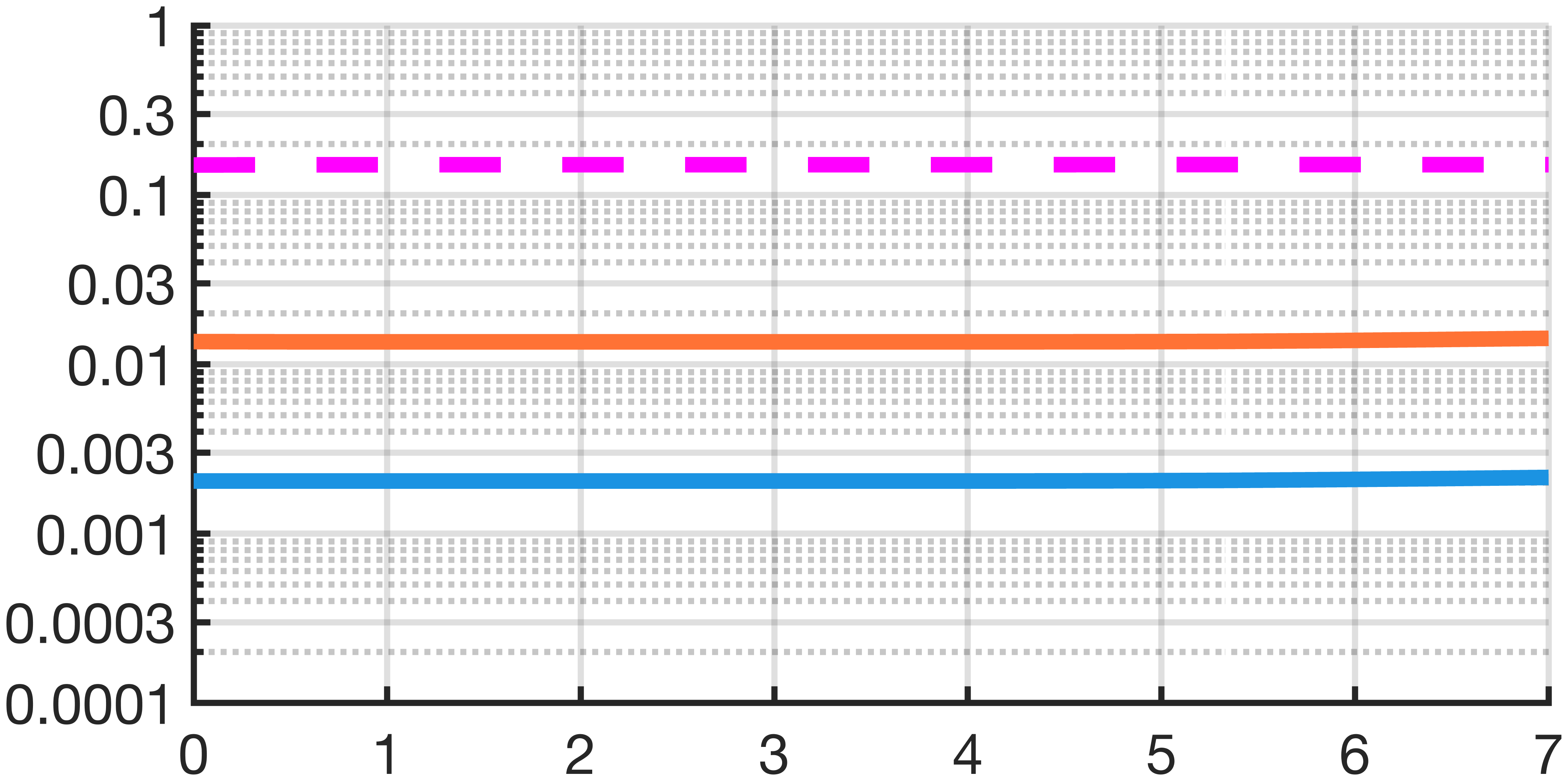} \\
WM-ROI \\
\end{minipage}
\begin{minipage}{5.5cm}
\centering
\includegraphics[height=2.5cm]{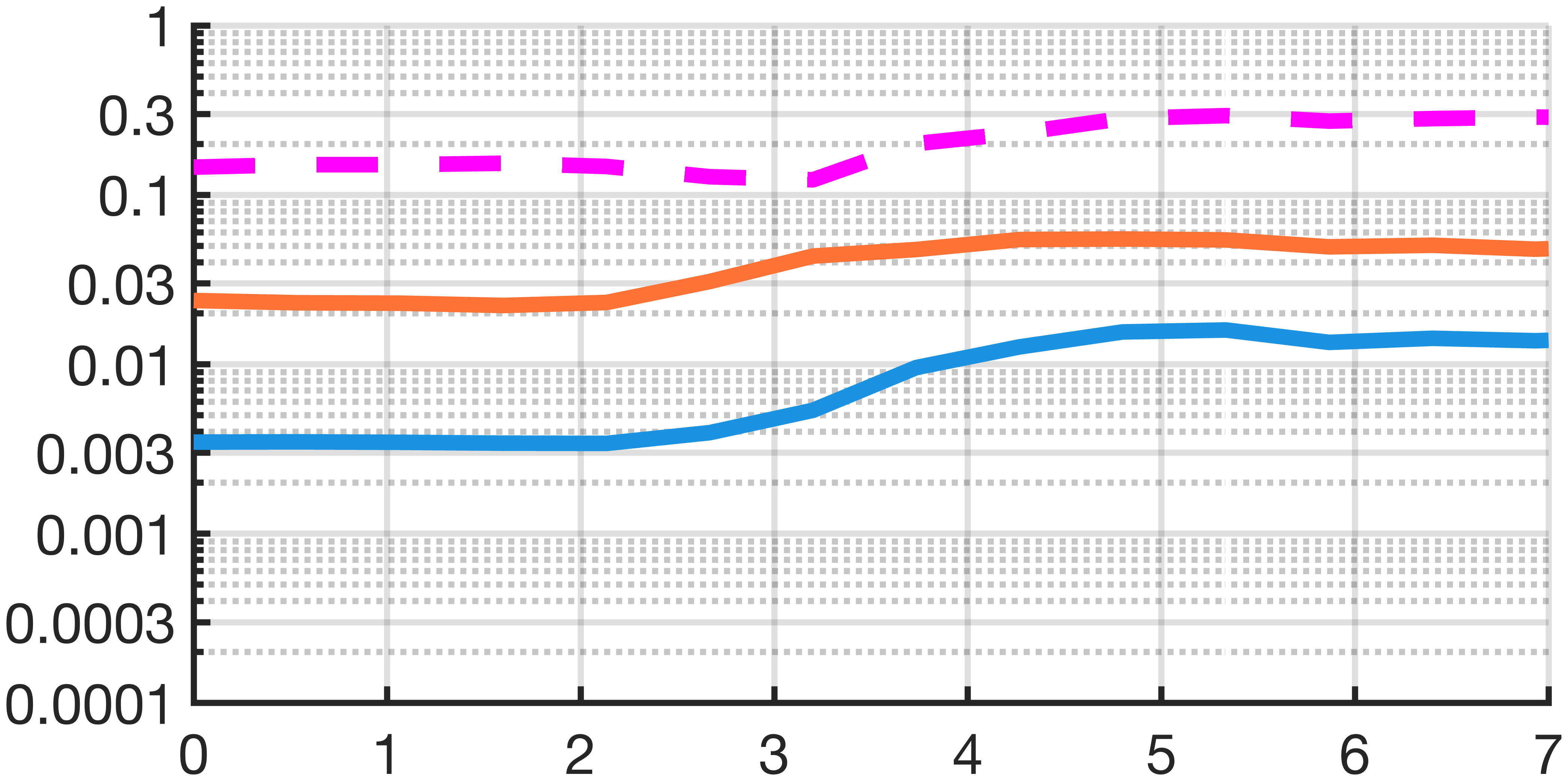} \\
GM-ROI \\
\end{minipage}
\begin{minipage}{5.5cm}
\centering
\includegraphics[height=2.5cm]{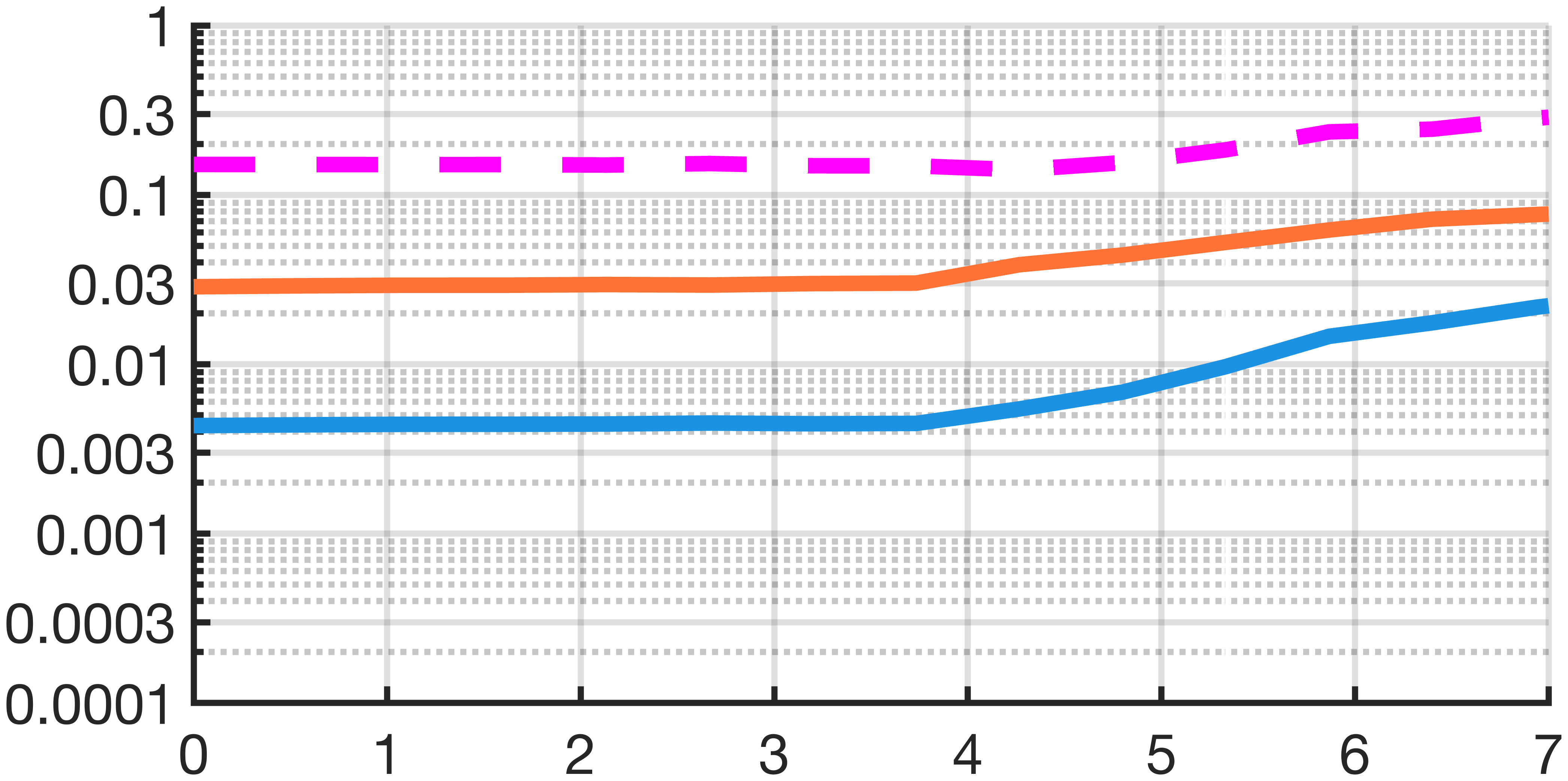} \\
BS-ROI \\
\end{minipage} \\
\vskip0.2cm
\begin{minipage}{5.5cm}
\centering
\includegraphics[height=2.5cm]{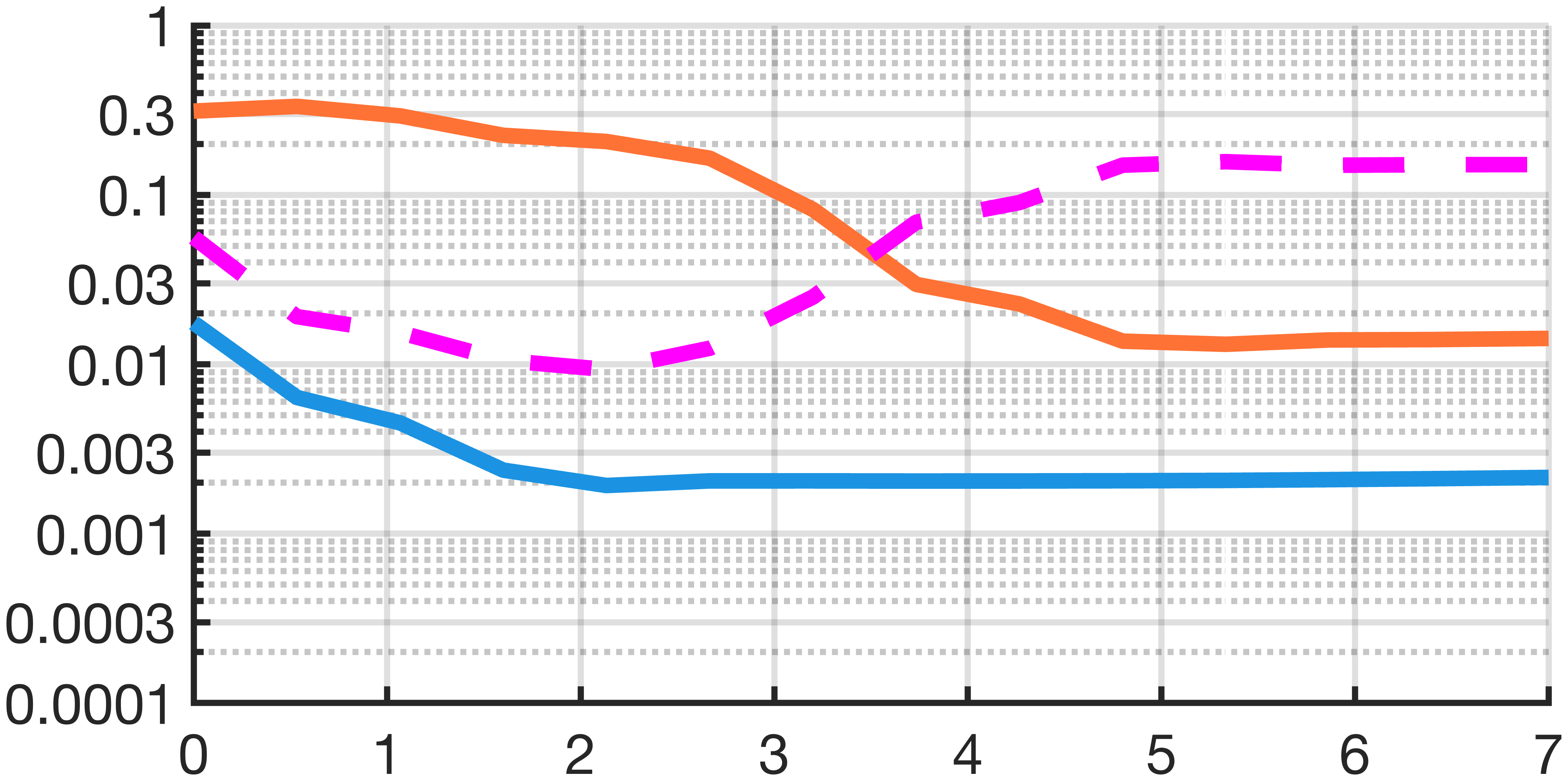} \\
WM-ROI \\
\end{minipage}
\begin{minipage}{5.5cm}
\centering
\includegraphics[height=2.5cm]{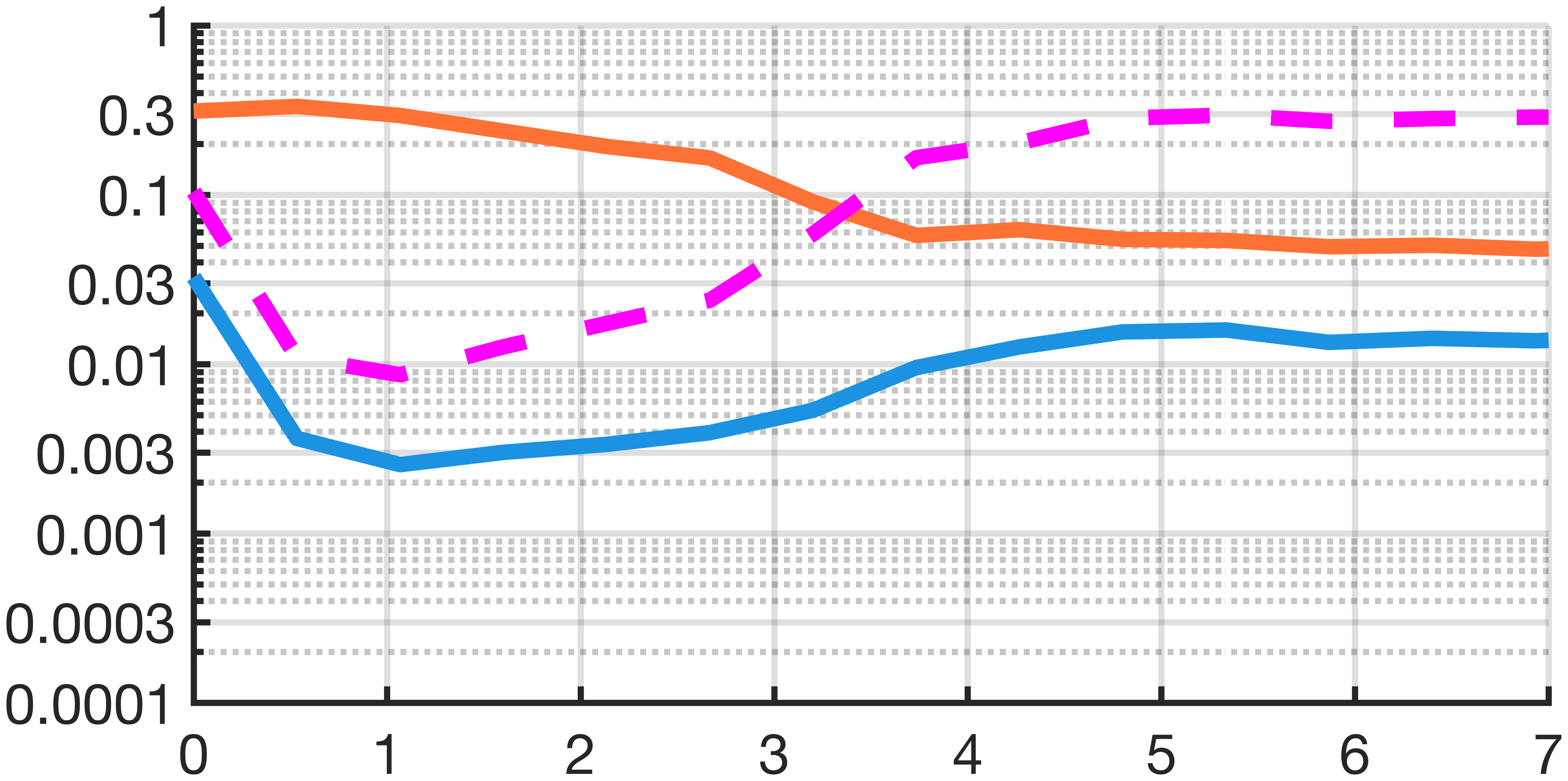} \\
GM-ROI \\
\end{minipage}
\begin{minipage}{5.5cm}
\centering
\includegraphics[height=2.5cm]{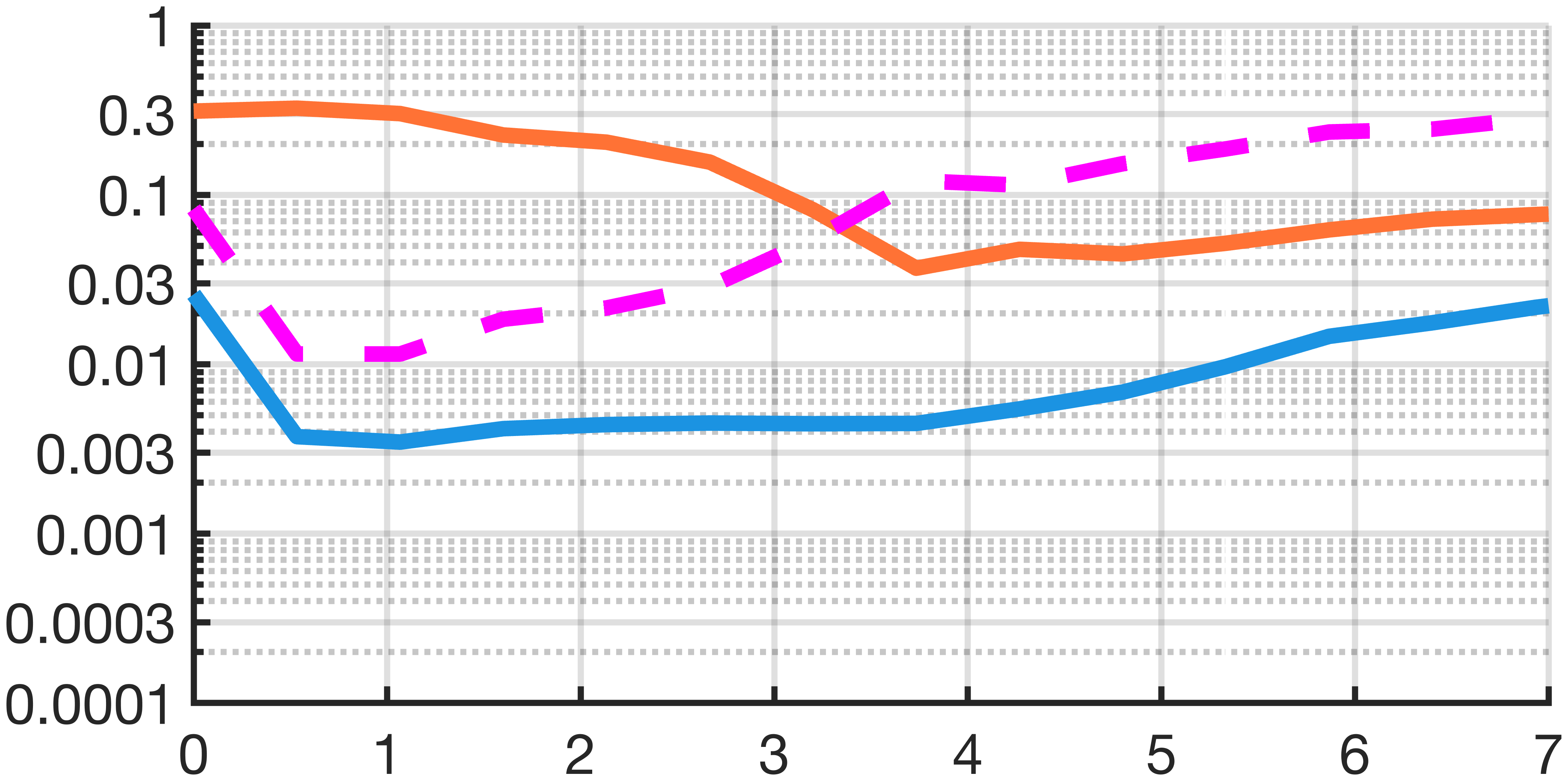} \\
BS-ROI \\
\end{minipage}
\end{footnotesize}
\caption{ The spatial profiles of the modeled total volumetric blood and deoxygenated blood volume (TBV and DBV) fractions within three ROIs at $t = 0.25$ s (top row) and $t = 5$ s (bottom row). The horizontal axis represents the distance from the center of the ROI, and the vertical axis shows the value of TBV (solid red line), DBV (solid blue line), and their ratio (DBV/TBV) (dashed magenta line). In each case, the neural activity is modelled as a point source in the center of the ROI and the value for a given time point has been obtained as integral (spatial) mean over the corresponding ROI. \label{fig:spatial_profiles}}
\end{figure*}

The spatial decay of TBV and DBV fractions, particularly in the WM-ROI, is slower compared to GM-ROI and BS-ROI, as illustrated in Figure  \ref{fig:spatial_profiles}. Notably, the WM-ROI exhibited the highest relative response, reflecting its higher metabolic demand and vascular density. In contrast, the GM-ROI, despite having a higher density of microvessels, demonstrated a suppressed relative response compared to the WM-ROI. The BS-ROI, with the highest microvessel density, exhibited the least response. The microvessel density influences the diffusion constant -- the higher density results in greater values and faster decay of the estimated blood volume fractions as the distance from the source of neural activity increases. This decay is more evident in the DBV fraction, which reduces faster than the TBV fraction. The weak response for GM-ROI and BS-ROI is likely due to the proximity of a major artery, which increases the baseline level of microcirculation further from the center of the ROI in those cases. This result is visible in Figure \ref{fig:spatial_profiles}, where the arterial flow significantly influences the DBV and TBV fractions profiles as early as $t = 0.25$ s, when the amplitude of the hemodynamic response remains negligibly small (Figure \ref{fig:oxy_vs_deoxy_temporal}). 

At $t=5$ s, the diameter of the perturbation in TBV fraction is between 8.0 and 8.5 mm, while in DBV fraction, it ranges from 3.4 to 4.0 mm (Table \ref{tab:physical_parameters}). The size of the perturbation is reflected by the decay curves in Figure \ref{fig:spatial_profiles} and seems to be influenced by variations in the diffusion coefficient, which demonstrate how variations in microvessel density influence the extent and rate of spatial diffusion. By comparing the relative amount of DBV and TBV fractions (Figure \ref{fig:spatial_profiles}), it can be observed that an approximate background level of 15 \% is maintained near the point of neural activity at the beginning of the hemodynamic response, as shown by the results at $t = 0.25$ s. Closer to the arteries, this relative amount grows to around 30 \%, or slightly higher, as seen in the results for the GM-ROI and BS-ROI (Figure \ref{fig:spatial_profiles}). However, in the absence of arteries, as indicated by the WM-ROI results, this value remains unchanged.

\begin{figure}[h!]
\centering
\begin{footnotesize}
\begin{minipage}{4.25cm}
\centering
\includegraphics[height=3cm]{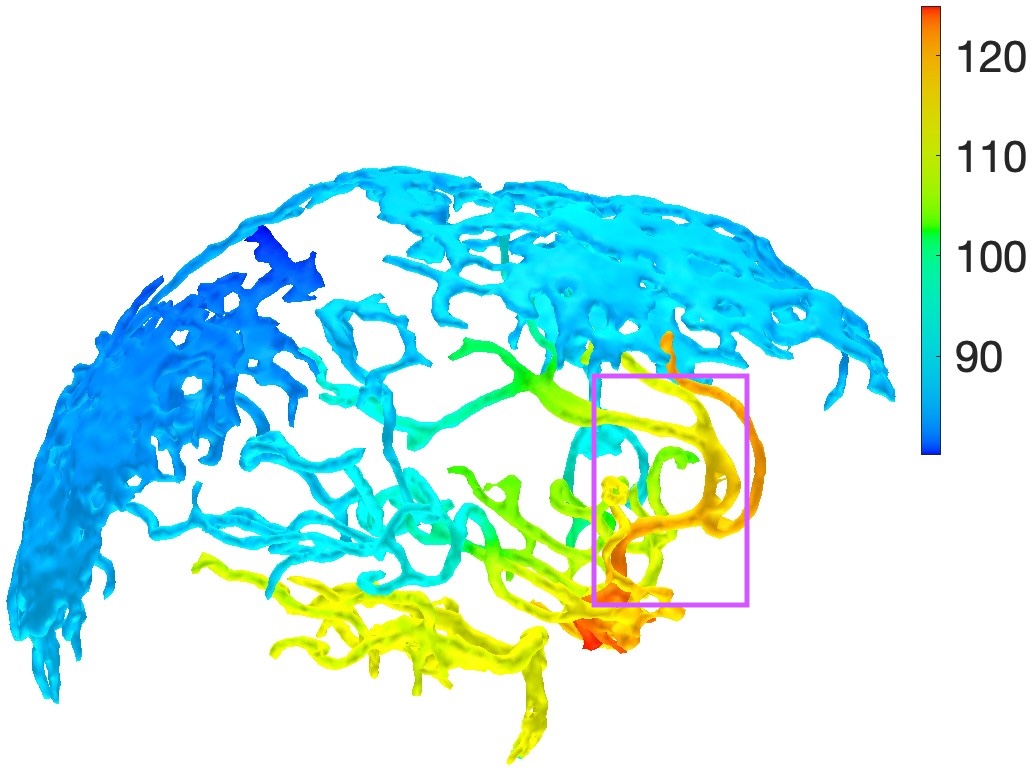} \\
Pressure 
\end{minipage}
\begin{minipage}{4.25cm}
\centering
\includegraphics[height=3cm]{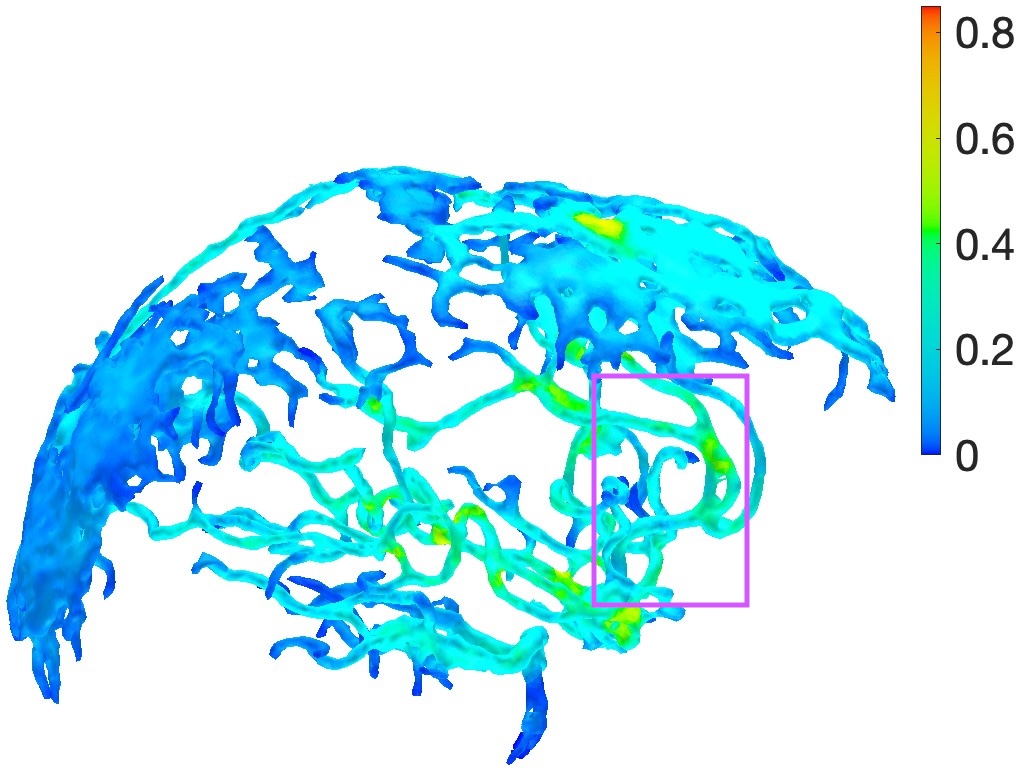} \\
Velocity 
\end{minipage} \\ \vskip0.1cm
\begin{minipage}{4.25cm}
\centering
\includegraphics[height=3cm]{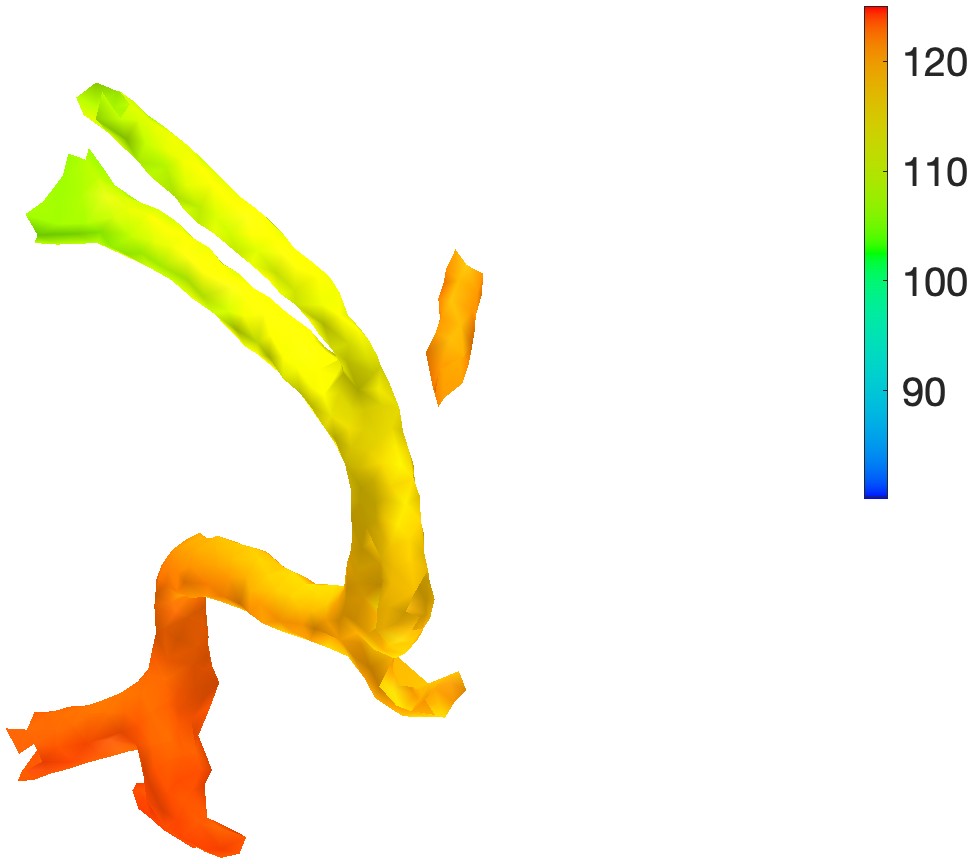} \\
Pressure in ACA
\end{minipage}
\begin{minipage}{4.25cm}
\centering
\includegraphics[height=3cm]{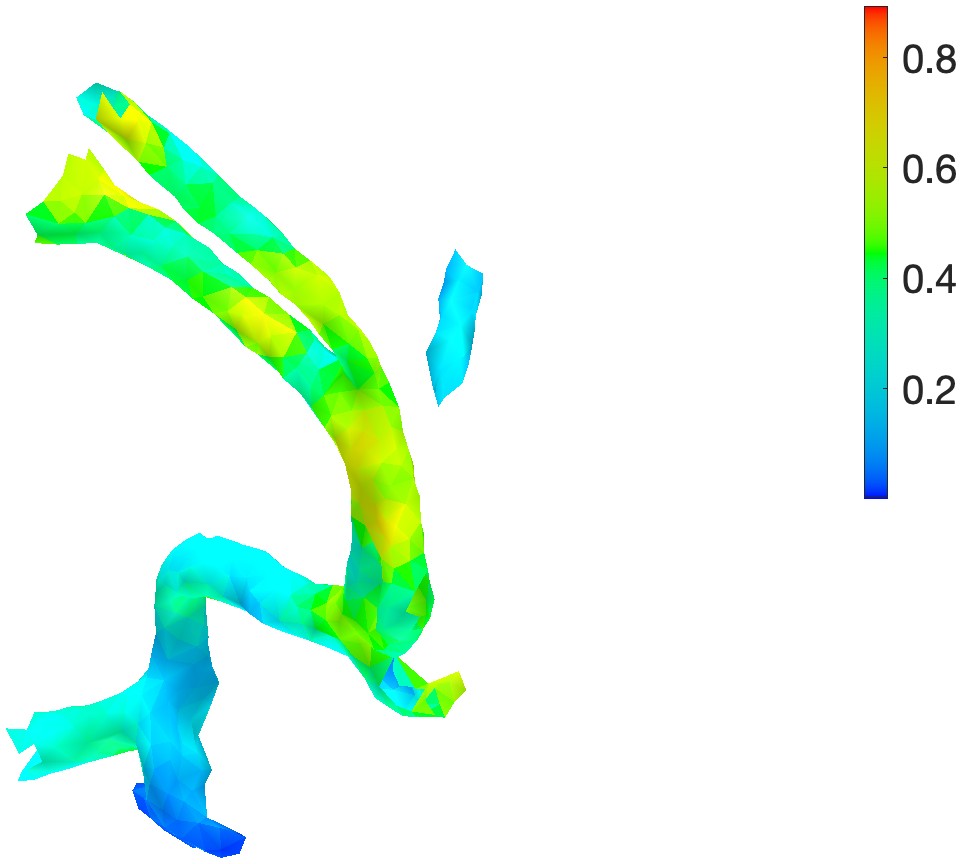} \\
Velocity in ACA
\end{minipage} \\ \vskip0.1cm
\begin{minipage}{4.25cm}
\centering
\includegraphics[height=3.3cm]{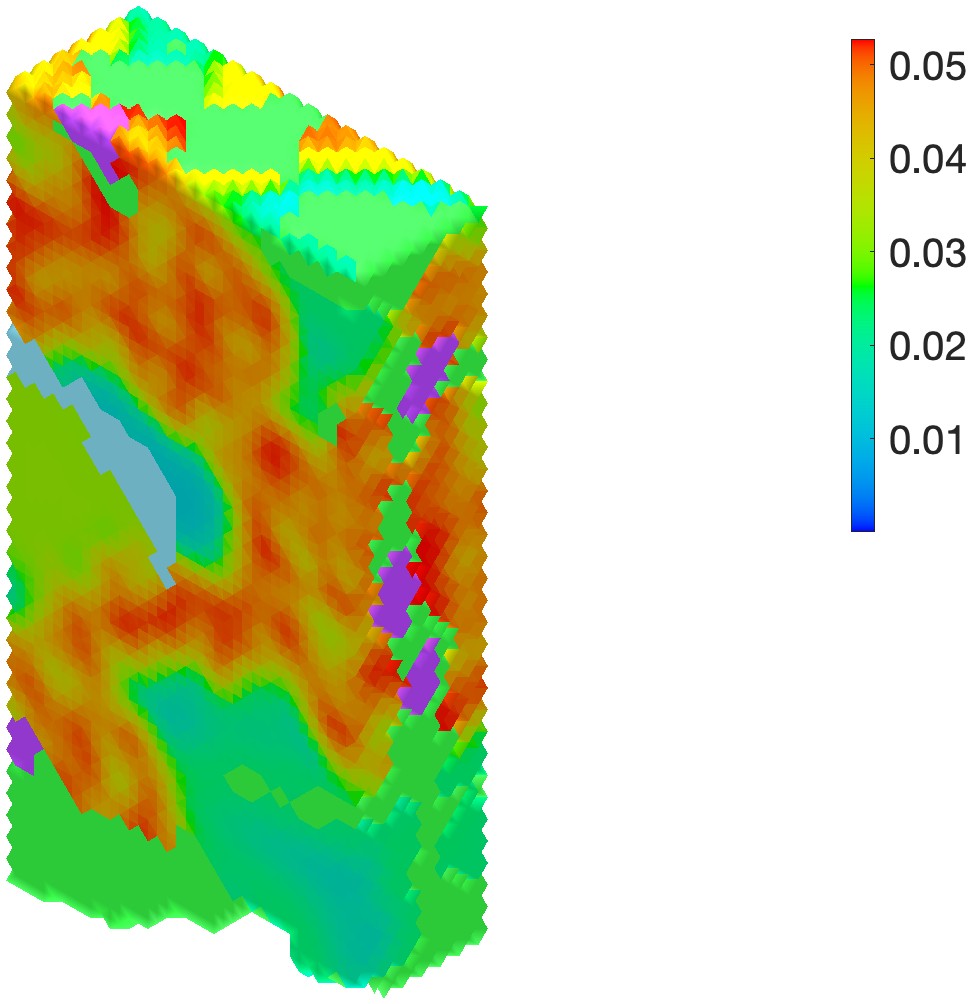} \\
TBV fraction around ACA
\end{minipage}
\begin{minipage}{4.25cm}
\centering
\includegraphics[height=3.3cm]{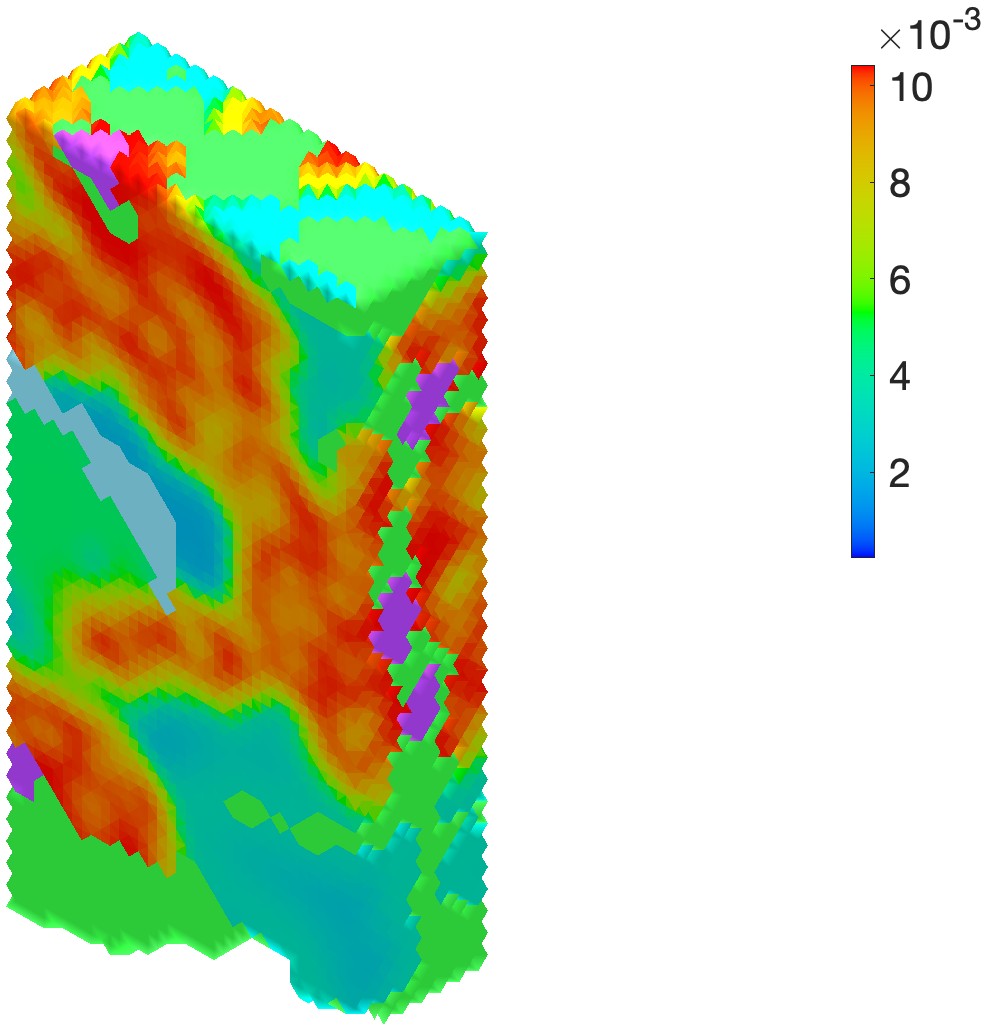} \\
DBV fraction around ACA
\end{minipage}
\\
\vskip0.1cm
\begin{minipage}{3cm}
\centering
\includegraphics[height=3.2cm]{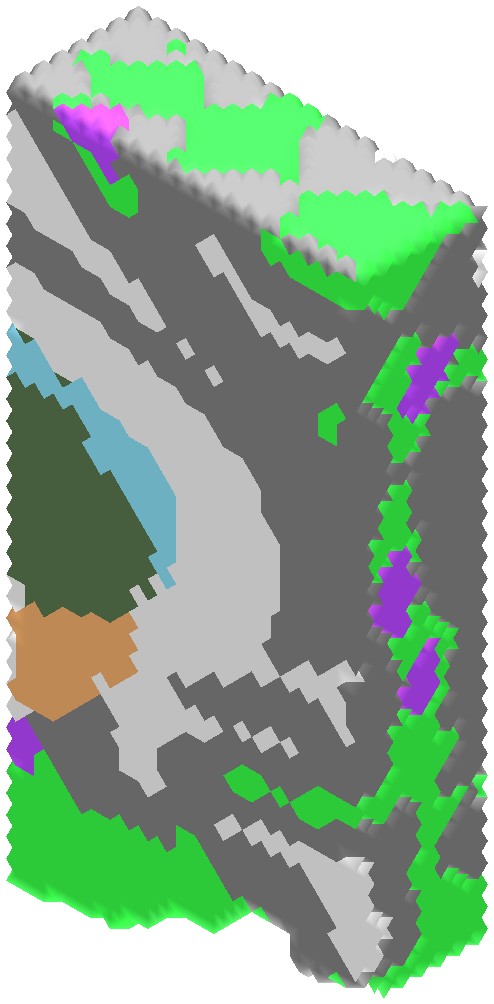} \\
Tissues
\end{minipage}
\begin{minipage}{5.5cm}
\centering
\includegraphics[width=4cm]{legend2.png}
\end{minipage}
\end{footnotesize}
\caption{Demonstration of the pressure (mmHg) and velocity (m/s) fields obtained as a solution PPE in the complete artery compartment (top row) and in the post-communicating, precallosal, and supracallosal parts of the anterior cerebral artery (ACA), for which the total and deoxygenated blood volume (TBV and DBV) fractions are also shown (third row). The corresponding tissue structure is shown in the bottom row. The box in the images of the top row shows the location of the ACA excerpt. \label{fig:aca}}
\end{figure}

An overall illustration of the pressure and velocity distributions obtained as a solution of PPE can be found in Figure \ref{fig:aca}. The relationship between those and the distributions of DBV and TBV is also shown, focusing in the post-communicating, precallosal and supracallosal parts of the anterior cerebral artery (ACA).

 % The white matter exhibited a significant response, which, while less pronounced than that in grey matter, was still notable given its lower vascular density. The brainstem, despite its high microvessel density, displayed the least response due to the influence of the neighboring artery. Coronal slices of TBV and DBV concentrations at $t = 5$ s  provide a visual representation of the spatial distribution of these quantities (Figures \ref{fig:VB_spatial}) and \ref{fig:dHb_spatial}). 

%%%%%%%%%%%%%%%%%%%%%%%%%%%
% \begin{figure*}[h!]
% \centering
% \begin{footnotesize}
% \begin{minipage}{5.5cm}
% \centering
% \includegraphics[height=2.5cm]{whitematter_area.png} \\
% WM-ROI \\
% \end{minipage}
% \begin{minipage}{5.5cm}
% \centering
% \includegraphics[height=2.5cm]{cortex_area.png} \\
% GM-ROI \\
% \end{minipage}
% \begin{minipage}{5.5cm}
% \centering
% \includegraphics[height=2.5cm]{brainstem_area.png} \\
% BS-ROI \\
% \end{minipage} 
% \end{footnotesize}
% \caption{The relative area of the set in which the relative difference to the local mean value of the volumetric OBV (red) or DBV (blue) concentration inside the ROI is greater than one. The horizontal axis shows the time $t$ (s). In each case, the neural activity is modelled as a point source in the center of the ROI.\label{fig:oxy_vs_deoxy_peak_area}}
% \end{figure*}
%%%%%%%%%%%%%%%%%%%%%%%%
%%%%%%%%%%%%%%%%%%%%%%%%
\begin{table*}[h!]
\caption{Relative volumetric blood (TBV) and deoxygenated blood volume (DBV) fractions and their ratio (DBV/TBV) for the tree different regions of interest (ROI) in the brain at $t=0$ and $t=5$ s. For the latter time point, the diameter of the perturbation caused by the hemodynamic response is given as the diameter in which the relative difference between the background and the distributions at $t=5$ s $t = 0.25$ s is greater than one inside the ROI.  \label{tab:bloodflow}}
\centering
\resizebox{0.6\linewidth}{!}{
\begin{tabular}{llrrrr}
\toprule
& & Background  & Volume &  Volume  & Perturbation \\
  & & volume  &  fraction & fraction & diameter  \\
    &  & fraction &  at $t = 0.25$ s & at $t = 5$ s & at $t = 5$ s \\
ROI   & Property &(relative)  &  (relative) & (relative) & (mm) \\
\midrule
{WM-ROI} & & & & \\
&TBV %(7)
& 0.0127  & 0.0138 & 0.0362 & 8.49 \\
&DBV %(8)
&        & 0.00208 & 0.00212 & 4.05\\
&DBV/TBV %(9)
&        & 0.151 & 0.0585 \\
{GM-ROI} & & & & \\
&TBV %(1) 
& 0.0200 & 0.0467  & 0.0729  & 8.42 \\
&DBV %(2)
&        & 0.0119 & 0.0119 & 4.02\\
&DBV/TBV %(3)
&        & 0.255   & 0.164 \\
{BS-ROI} & & & & \\
&TBV %(4) 
& 0.0274 & 0.0530 &  0.0733 & 7.96   \\
&DBV %(5)
&        & 0.0112 &0.0112 & 3.35 \\
&DBV/TBV %(6)     
&        & 0.211 & 0.153 \\
\bottomrule
\end{tabular}
} % end resizebox
\end{table*}
%%%%%%%%%%%%%%%%%%%%%
%%%%%%%%%%%%%%%%%%%%%

%%%%%%%%%%%%%%%%%%
%%%%%%%%%%%%%%%%%%
%%%%%%%%%%%%%%%%%%
%

\section{Discussion}
\label{sec: Discussion}
The results of this study demonstrate that our coupled model of the eccess total blood volume (TBV) and deoxygenated blood volume (DBV) fractions, based on Fick’s law \cite{berg2020modelling, arciero2017mathematical, reichold2009vascular}, effectively captures the complex interplay between blood flow, oxygen transport, and neural activity. By integrating volumetric sources (neural activity), boundary sources (arterial inflow), and venous circulation sinks, the model goes beyond existing convolution models \cite{aquino2014spatiotemporal,lacy2022cortical,pang2017effects,pang2016response,friston2000nonlinear, shmuel2007spatio}, offering a more realistic representation of the spatiotemporal dynamics in brain microcirculation.

As our study primarily focuses on the mathematical modelling of cerebral hemodynamics, the validation of the model was conducted through comparison with existing empirical data from previous studies rather than through direct clinical validation. Our computational and modelling results for the hemodynamic response are consistent with biomedical data \cite{ wolf_2002, Venclove_2015}. Specifically, our simulations show a maximum increase in oxygenated blood volume (OBV) fraction at approximately $t = 5$ s, which closely aligns with the observed data from \cite{Venclove_2015}. After reaching this peak, OBV fraction levels decrease while DBV fraction levels increase, both gradually returning to baseline levels. These simulations also reveal differences in the hemodynamic responses across the regions of interest (ROIs), indicating spatial variations in the vascular and metabolic properties of different brain regions.

In this context, the diffusion constant plays a critical role in determining the rate at which TBV and DBV decay across the domain. Higher diffusion constants lead to faster and more pronounced decay. Numerical results indicate that the decay from a point-like neural source to the background value occurs at approximately 4.0 – 4.3 mm for TBV and 1.7 – 2.0 mm for DBV. Subtle differences are observable across ROIs, with WM-ROI, GM-ROI, and BS-ROI showing varying decay profiles. Our coarse approximation of the diffusion constant was derived using laminar Hagen-Poiseuille flow, consistent with previous studies that have modeled TBV decay using convolution-based approaches \cite{aquino2014spatiotemporal, lacy2022cortical, pang2016response, pang2017effects, pang2018biophysically, aquino2014deconvolution, caro2012mechanics}. Additionally, findings from \cite{shmuel2007spatio} on the spatiotemporal point-spread function of the BOLD MRI signal, particularly in the context of 7T MRI, show similar spatial extents to the DBV values observed in our model.

Based on these observations, our method provides a reasonable approximation not only for the decay of TBV near vessel boundary sources \cite{samavaki2023pressure}, but also for both TBV and DBV in proximity to point-like sources of hemodynamic response due to neural activity. The model assumes that regions with increased blood supply correspond to areas of higher metabolic activity, which is reflected in the growth of TBV and DBV profiles toward both point-like and boundary sources. In particular, the relative amount of DBV can increase from the background assumption of approximately 15\% to around 30\%, or slightly higher, a result consistent with experimental data on oxygen extraction across different brain regions \cite{ito2023oxygen}.

The use of 7T MRI data was essential for ensuring the accuracy of our vessel reconstruction and model application. The high-resolution imaging allowed us to distinguish brain compartments with precision, particularly in differentiating small vessels from surrounding tissue, a task that lower-field MRI scans may struggle to accomplish \cite{svanera2021cerebrum, choi2020cerebral}. This enabled the model to achieve greater accuracy in simulating cerebral hemodynamic responses, consistent with empirical findings from related biomedical studies \cite{wolf_2002, Venclove_2015}.

\subsection{Limitations of the Study}

Our study is fundamentally a mathematical exploration of diffusion-driven processes in CBF. While the model's predictions align with empirical findings, further validation is necessary in a clinical context, particularly through comparisons with high-resolution perfusion imaging or other direct measurements of CBF and metabolism. Incorporating additional empirical data and performing validation studies across various neurological conditions are essential steps toward advancing the model from a theoretical construct to a clinically applicable tool.

The primary limitation of our approach lies in the spatial modelling accuracy of the diffusion-driven framework. While our method successfully captures essential aspects of the spatiotemporal dynamics of CBF and oxygenation, its spatial precision does not yet meet the standards required for clinical applications, such as those in perfusion imaging studies, see, e.g., \cite{zhang2022pilot,chen2018evaluation}. In clinical contexts, highly accurate modelling of cerebral perfusion is crucial for diagnostics and treatment planning, particularly for conditions like stroke and tumors.  Thus, further refinements, including improved parameterization and validation against clinical perfusion imaging data, are needed to bridge the gap between mathematical modeling and clinical practice.

As for the limitations of the mathematical model lies in the assumptions underlying the application of Fick's law, particularly the assumption of a dense microvascular network. This approach omits global convection and relies primarily on diffusion as the main transport mechanism. Although this simplifies the model, it may not fully capture regions of the brain where convection plays a more significant role, such as in cases of abnormal pathology or areas with high blood flow.

Additionally, the approximation of the diffusion constant using average parameters from the Hagen-Poiseuille flow formula introduces a degree of uncertainty into the model. This approach simplifies the complicated reality of blood flow by assuming homogeneity, which may overlook minor but significant differences in tissue structure, vessel size, and distribution. Such variations can significantly affect local hemodynamics, influencing how blood and oxygen are delivered to specific regions of the brain \cite{caro2012mechanics}. 

Furthermore, our model does not account for potential variability in the oxygen extraction fraction and the cerebral metabolic rate of oxygen. These factors can exhibit significant differences in their distribution across various brain areas and physiological situations. By neglecting these variations, the model may produce discrepancies when compared with experimental data at finer spatial resolutions. Incorporating region-specific oxygen extraction fraction and cerebral metabolic rate of oxygen values would be essential to improve the model’s accuracy, especially in scenarios where localized metabolic activity or oxygen dynamics play a critical role \cite{ito2023oxygen}.

%%%%%%%%%%%%%%%%%%%%%%
%%%%%%%%%%%%%%%%%%%%%%
%%%%%%%%%%%%%%%%%%%%%%
\subsection{Broader Implications and Impact}

Our findings contribute to the expanding body of research on CBF modelling by by highlighting the benefits of integrating diffusion-driven processes with advanced imaging techniques for more accurate cerebral hemodynamic modelling. This approach holds promise for mathematical, experimental, and clinical research, particularly in the context of BOLD fMRI and NIRS imaging  \cite{ozaki2021near, puckett2016spatiotemporal,kinou2013differential, van2020dynamics}. It offers valuable insights into conditions where blood flow and oxygen transport are compromised, such as stroke, traumatic brain injury, and neurodegenerative diseases.

Moreover, incorporating a realistic head model derived from 7T MRI data \cite{svanera2021cerebrum} offers significant potential for advancing the accuracy of brain dynamics simulations. This, in turn, could lead to enhanced diagnostic and therapeutic tools for medical professionals. For example, more precise models of blood oxygenation may support the development of targeted interventions for conditions that impact CBF.

%%%%%%%%%%%%%%%%%%%%%%%%
%%%%%%%%%%%%%%%%%%%%%%%%
%%%%%%%%%%%%%%%%%%%%%%%%

\subsection{Specific Takeaways and Future Research Directions}

This study makes two key contributions: (1) the development of an advanced model that incorporates diffusion-driven transport for both TBV and DBV, and (2) the validation of this model using high-resolution data, demonstrating its significance for neuroimaging applications and physiological research.

Future research should focus on addressing the current limitations of the model. A promising avenue for refinement is the incorporation of convection in areas where it may play a more significant role. Additionally, efforts could be made to improve the estimation of the diffusion constant by incorporating data on local tissue properties and microvessel distributions, potentially through in vivo measurements or more detailed post-mortem analyses.

Another important direction for future research is the development of more sophisticated models of oxygen extraction and consumption, moving beyond the bulk approximation of the balloon model toward a more detailed metabolic framework \cite{buxton2004modeling, miraucourt2017blood}. Incorporating dynamic changes in blood viscosity, particularly in response to fluctuations in blood oxygenation, could also lead to more accurate predictions of hemodynamic responses \cite{valant2016influence, KANARIS_2012}.

%%%%%%%%%%%%%%%%%%%%%%
%%%%%%%%%%%%%%%%%%%%%%
%%%%%%%%%%%%%%%%%%%%%%
\subsection{Conclusion and Call to Action}

From a mathematical modelling perspective, our research strongly supports the use of diffusion-driven approaches for simulating CBF. The model employed in this study provides a more detailed mathematical representation of spatiotemporal dynamics within brain microcirculation, particularly in areas where convection is negligible. Specifically, the diffusion approximation based on the Hagen-Poiseuille model was proposed to account for the point spread of TBV and DBV responses observed in experimental studies.

Our approach provides a new avenue for investigating the  mechanisms of blood flow and oxygen transport in the brain, offering a valuable tool for understanding the physiological processes underlying neurovascular coupling. While our model aligns with empirical findings and highlights the potential of diffusion-driven approaches, it is important to acknowledge that its accuracy is not yet sufficient for direct clinical applications. The coarse approximation of certain parameters, such as the diffusion constant, and the assumptions within the model limit its current clinical utility. Nevertheless, this work establishes a foundation for future refinements that could improve its applicability in medical diagnostics and treatment planning.

We encourage the integration of diffusion-driven models in future studies of CBF, particularly in contexts where spatial resolution and accuracy are critical.  With continued refinement, including the incorporation of more precise parameter values and the consideration of convection where  appropriate, this modelling approach could evolve into a powerful tool for investigating a broader range of neurological conditions. Although not yet ready for clinical application, these models show great potential for advancing our understanding of cerebral hemodynamics and for contributing to the development of improved diagnostic and therapeutic strategies.

%%%%%%%%%%%%%%%%%%%%%%%%%%
%%%%%%%%%%%%%%%%%%%%%%%%%%
%%%%%%%%%%%%%%%%%%%%%%%%%%

\appendix
%{\color{black}
%\section{Length Density of Arterioles}
%\label{app:arteriole_length_density}

%%%%%%%%%%%%%%%%%
%%%%%%%%%%%%%%%%%
%%%%%%%%%%%%%%%%%
\section*{Fundings and Acknowledgements}
The work of Maryam Samavaki and Sampsa Pursiainen was supported by the Research Council of Finland's Flagship of Advanced Mathematics for Sensing Imaging and Modelling 2024-2031 (grant 359185) and Centre of Excellence in Inverse Modelling and Imaging 2018-2025 (grant 353089). Santtu Söderholm and Maryam Samavaki are supported by the ERA PerMed project PerEpi (PERsonalized diagnosis and treatment for refractory focal paediatric and adult EPIlepsy), RCF's decision 344712; Arash Zarrin Nia has been funded by a scholarship from the K. N. Toosi University of Technology. We thank Prof. Carsten H. Wolters, University of Münster, Münster, Germany, for the fruitful discussions and support in mathematical modelling and are grateful to DAAD (German Academic Exchange Service) and RCF for supporting our travels to Münster (RCF decision 354976).
%The work of  is supported by the Academy of Finland Centre of Excellence (CoE) in Inverse modelling and Imaging 2018–2025 (decision 336792) and project 336151; Santtu Söderholm is funded by the Research Council of Finland (previously Academy of Finland) project \#344712 (Era PerMed/PerEpi); and  has been funded by a scholarship from the K. N. Toosi University of Technology.
%DAAD 354976,Center of Exellent 353089, Flagship 359185, PerEpi 344712

%%%%%%%%%%%%%%%%%%%%%%%%
%%%%%%%%%%%%%%%%%%%%%%%%
%%%%%%%%%%%%%%%%%%%%%%%

\section*{Conflicts of Interest}
\label{sec:conflict}
The authors confirm that the research utilized in this study was entirely independent, open, and academic.
They have neither financial nor non-financial relationships, affiliations, knowledge, or beliefs in the
subject matter or materials included in this manuscript, nor do they have any involvement with or
affiliation with any organization or personal relationship.

%%%%%%%%%%%%%%%%%%%%%
%%%%%%%%%%%%%%%%%%%%%%
%%%%%%%%%%%%%%%%%%%%%

\section*{Ethical Considerations}
\label{sec:ethics}
This study only involves mathematical modelling and numerical simulations, which utilize an openly accessible MRI dataset from \href{https://openneuro.org/}{OpenNeuro}~\cite{our-7T-dataset} as input. Therefore, the conductors of this study did not gather or process any sensitive data at any phase of the study, and hence no ethical guidelines could have been violated in the process. No human or animal subjects were involved in this research.

%%%%%%%%%%%%%%%%%%%
%%%%%%%%%%%%%%%%%%%
%%%%%%%%%%%%%%%%%%%

\bibliographystyle{elsarticle-num} 
\bibliography{cas-refs}

%% else use the following coding to input the bibitems directly in the
%% TeX file.

% \begin{thebibliography}{00}

% %% \bibitem{label}
% %% Text of bibliographic item

% \bibitem{}

% \end{thebibliography}
\end{document}